\documentclass{amsart}
\usepackage[active]{srcltx}
\usepackage{amsmath,amsthm,amscd,amssymb}
\usepackage{latexsym}
\usepackage{floatflt}
\usepackage{wrapfig}
\usepackage[final]{graphicx}
\usepackage{pst-ghsb}
\usepackage{epic, eepic}
\def\squarebox#1{\hbox to #1{\hfill\vbox to #1{\vfill}}}

\newcommand{\cz}{{\mathbb C}}
\newcommand{\nz}{{\mathbb N}}
\newcommand{\qz}{{\mathbb Q}}
\newcommand{\rz}{{\mathbb R}}
\newcommand{\zz}{{\mathbb Z}}

\def\tr{\mathop{\mathrm{tr}} \nolimits} 
\def\qed{\hbox {\hskip 1pt \vrule width 4pt height 6pt depth 1.5pt
        \hskip 1pt}}

\def\rhs{right hand side}

\def\Re{{\rm Re\,}}
\def\Im{{\rm Im\,}}
\newtheorem{theorem}{Theorem}

\newtheorem{lemma}{Lemma}

\newtheorem{proposition}{Proposition}

\title{Resonance spectrum  for one-dimensional layered media}

\author{Alexei Iantchenko}
\address{Malm\"o University\\
School of Technology and Society\\
SE-205 06 Malm\"o\\
Sweden } \email{ ai@ts.mah.se}
\date{\today}

\begin{document}

\bibliographystyle{plain}



\begin{abstract}
We consider the ``weighted'' operator $P_k=-\partial_x
a(x)\partial_x$ on the line with a step-like coefficient which
appears when propagation of waves thorough a finite slab of a
periodic medium is studied. The medium is transparent at certain
resonant frequencies which are related to the complex resonance
spectrum of $P_k.$

If the coefficient is periodic on a finite interval (locally
periodic) with $k$ identical cells then the resonance spectrum of
$P_k$ has band structure. In the present paper we study a transition
to semi-infinite medium by taking the limit $k\rightarrow \infty.$
 The bands of resonances in the complex lower half plane are
 localized
 below the band spectrum of the corresponding periodic problem
($k=\infty$) with  $k-1$ or $k$ resonances in each band. We prove
that as $k\rightarrow \infty$ the resonance spectrum converges to
the real axis.
\end{abstract}

\pagestyle{myheadings} \maketitle

\section{Introduction}

In the present paper we consider operator $P_k=-\partial_x
a_k(x)\partial_x$ on the line with step-like coefficient  $a_k$
which is periodic on a finite interval defined as follows:
\begin{equation}\label{pk}
a_k(x)=a(x),\,\,\mbox{for}\,\,x\in
[0,k];\,\,a_k(x)=\frac{1}{b_1^2},\,\,\mbox{for}\,\,x\not\in [0,k]
,\end{equation} where $a(x)$ is $1-$periodic function equal to
\begin{equation}\label{a0}
a_0 (x)=\left\{\begin{array}{lr}
                   b_2^{-2} & \mbox{for}\,\,x\in [0, x_2) \\
                   b_1^{-2} & \mbox{for}\,\, x\in [x_2,1)
                 \end{array}\right.\end{equation}
                 for $x\in [0,1).$ Here $b_{1,2}>0$ and $0<x_2<1.$
Equation
$$ P_k\psi =-\partial_x a_k(x)\partial_x\psi(x)=\lambda^2\psi$$
 appears when the propagation of waves through a finite
slab of a periodic medium is studied. Such systems are also called
finite or locally periodic media (for revue see
\cite{GriffithsSteinke2001}).

When $k$ is large then the properties of medium is close to an
infinite periodic problem in a sense that we are going to discuss in
the present paper.

We denote  $P=-\partial_x a(x)\partial_x$ the pure periodic
operator, where
 $a(x)$ is
$1-$periodic function equal to $a_0$ for $x\in [0,1)$ as in
(\ref{a0}). Then the Floquet theory shows  the existence of a pair
of the quasi-periodic solutions $\psi_\pm$ of the equation
$-\partial_x a(x)\partial_x\psi_\pm =\lambda^2\psi_\pm,$
$$\psi^\pm(\lambda,x+1)=e^{\pm i\theta}\psi^{\pm}(\lambda,x),$$
such that $\psi^\pm\in L^2(\rz_\pm)$ for $\Im\lambda >0.$ Here
$\theta=\theta(\lambda)$ is the Bloch phase. We denote
\begin{equation}\label{lyaponov}
F(\lambda)= \frac{\rho+1}{2}\cos\{\lambda(x_2b_2
+(1-x_2)b_1)\}-\frac{\rho-1}{2}\cos\{\lambda(x_2b_2-(1-x_2)b_1)\}
\end{equation} the Lyaponov function for $P$ (see Section
\ref{ss-Lyaponov}). Here
$$\rho=\frac{b_1^2+b_2^2}{2b_1b_2}.$$ The spectrum of the operator $P$ has band
structure with allowed zones defined as follows:
\begin{equation}\label{bandspectrum}
\lambda\in\sigma(P)\,\,\Leftrightarrow\,\,|F(\lambda)|^2
<1,\,\,\lambda\in\rz
\end{equation}
  (see \cite{Banica2003} and
Section \ref{s-banica}). The band edges are given by solutions of $
F(\lambda)=\pm 1.$

The relation between the Bloch phase $\theta$ and the spectral
parameter $\lambda$ is called dispersion relation: $$\cos\theta
(\lambda)=F(\lambda).$$

Since the coefficient $a(x)$ is constant equal to $1/b_1^2$ outside
a finite region, we are here concerned with a scattering problem.

 We shall denote the reflection
and transmission coefficients for the operator $P_k$  by $r_k$ and
$t_k,$ respectively:
$$P_k\psi =\lambda^2\psi,\,\,\psi=e^{i\lambda b_1x} +r_k e^{-i\lambda
b_1 x},\,\,x <0;\,\,\psi=t_ke^{i\lambda b_1x} ,\,\,x >k.$$

Following the ideas in \cite{MolchanovVainberg2006} we consider a
transition to semi-infinite periodic materials by taking the limit
$k\rightarrow\infty$ of the reflection coefficient $r_k$ for $P_k.$

The limiting operator
$$P_\infty\psi =-\partial_x a_\infty(x)\partial_x\psi (x)$$
corresponds to the case of such a long slab that it can be
considered as half infinite.

In the case of the operator $P_\infty,$ the solution $\psi$ of the
scattering problem is defined as the solution of the equation
$P_\infty\psi =\lambda^2\psi,$ such that
\begin{equation}\label{(190)}\psi=e^{i\lambda b_1x} +r e^{-i\lambda b_1 x},\,\,x
<0;\,\,\psi= c \psi^+(\lambda, x) ,\,\,x
>0,\end{equation} with some $r=r(\lambda),$ $c=c(\lambda).$

As in \cite{MolchanovVainberg2006} we have that
 the reflection coefficients $r_k(\lambda)$ and
$r(\lambda)$ are analytic in the upper half plane
$\cz_+=\{\lambda:\,\,\Im\lambda >0\}$ and continuous in
$\overline{\cz}_+,$ and $r_k(\lambda)\rightarrow r(\lambda)$ when
$k\rightarrow\infty$ and $\lambda\in\cz_+.$ When $\lambda$ is real,
$r_k(\lambda)$ converges to $r(\lambda)$ in the weak sense (see
Theorem \ref{tMV}, Section \ref{ss-transition}).

Numerical calculation  shows  that in each allowed zone of
$\sigma(P)$ there are in general  $k-1$  frequencies $\lambda_j$
where the transmission probability is one:
$|t_k(\lambda_j)|^2=1,\,\,j=1,\ldots,k-1,$ and the medium is
perfectly transparent: $|r_k(\lambda_j)|^2=0.$ There exist an
additional frequency $\lambda_0$ when the medium consisting of only
one unit cell is transparent and then $|t_n(\lambda_0)|^2=1$ for all
$n=1,2,\ldots,k.$
 The pics in the transmission probability are
 related to the  complex resonances close to the real
axis.

We make the following definition.

The operator $P_k$ defined from $\{ u\in H^1(\rz),\,\,a\partial_x
u\in H^1(\rz)\}$ to $L^2(\rz)$ is self-adjoint.
For $\Im\lambda > 0,$  we call $R_k(\lambda) v=(-\partial_x
a(x)\partial_x -\lambda^{2})^{-1} v$ the resolvent of $P_k.$ For any
$k=1,2,\ldots,$ the operator-valued function
$$R_k(\lambda):\,\,L^2_{\rm comp} (\rz)\mapsto L_{\rm loc}^2(\rz)$$
can be continued to the lower complex half-plane $\cz_-$ as a
meromorphic function of $\lambda\in\cz$ and it has no poles for
$\Im\lambda \geq -\epsilon_k,\,\,\lambda\neq 0, $ with
$\epsilon_k>0$ positive constant dependent on $k$ (see Section
\ref{s-banica}).

The poles of the $R_k(\lambda)$ in $\cz_-$ are called {\em
resonances} or {\em scattering poles}. We denote the set of
resonances ${\rm Res}\,(P_k).$

 Using the explicit construction of the resolvent in
 \cite{Banica2003}  the poles are calculated numerically.
Some examples are presented in Section \ref{s-numerics}, figures
(\ref{Fig1}), (\ref{Fig2}) and (\ref{Fig3}). We summarize the
properties of ${\rm Res}\,(P_k)$ in the following Theorem.
\begin{theorem}\label{main} We consider the finite periodic $P_k$ and
periodic $P$ operators generated by the same unit cell given in
(\ref{a0}).
 Let ${\rm Res}\,(P_k)\subset \cz_-$ denote the resonance spectrum for the finite periodic system with   $k=2,3\ldots$ identical cells and
$\sigma(P)\subset\rz$ denote the band spectrum for $P$ given by
(\ref{bandspectrum}).

   The resonance spectrum for the finitely periodic system has band structure related to the bands of the real spectrum for the pure periodic problem as
   follows:\\
{\bf 1)} The resonance spectrum of $P_k$ has band structure.
Resonances are localized below the bands of the real spectrum of
$P:$
$$\lambda\in{\rm
Res}\,(P_k)\,\,\,\,\Rightarrow\,\,\,\,\Re(\lambda)\,\,\mbox{satisfies}\,\,
(\ref{bandspectrum})\,\,\,\,\Leftrightarrow\,\,\,\,
\Re(\lambda)\in\sigma(P_\infty).$$  Each resonance band of $P_k$
consists of  $k-1$ resonances $\lambda_1,\ldots,\lambda_k-1$ and
eventually an additional resonance with real part
$\Re\lambda=\lambda_0=\pi m/x_2b_2,$ $m\in\zz,$ such that
$|t_1(\lambda_0)|^2=1,$ the one-cell medium is ``perfectly
transparent'' at frequency $\lambda_0.$\\
{\bf 2)} If the condition \begin{equation}\label{cond}
b_2x_2=b_1(1-x_2)\,\,\Leftrightarrow\,\,\frac{b_1}{x_2}=\frac{b_2}{1-x_2
}.\end{equation} is satisfied  then $\lambda_0=\pi m/x_2b_2,$
$m\in\zz,$ is the degenerate band edge (two bands has common edge at
$\lambda_0$). The resonance spectrum ${\rm Res}\,(P_k),$
$k=1,2,3,\ldots,$ is periodic with the period $T=\frac{\pi}{b_2
x_2}.$\\
{\bf 3)} As $k\rightarrow\infty$ then the resonance spectrum of
$P_k$ approaches
 the real axis.

\end{theorem}

In the present paper we motivate these numerical results.

The band structure of the  resonance spectrum for a finitely
periodic system  and its relation to the band spectrum of the
correspondent periodic problem is well-known in physical literature
(see \cite{Cohen-Tannoudjietal1977}).

We say that ${\rm Res}\,(P_k)$ is periodic if there exists $T>0,$
period, such that
$${\rm Res}\,(P_k)\cap ([q+Tn,p+Tn]-i\rz)={\rm Res}\,(P_k)\cap
([q+Tm,p+Tm]-i\rz) $$ for any $q <p$ and $n,m\in\zz.$ This property
follows directly from the equations defining the resonances  in
Section (\ref{ss-period}) if  condition (\ref{cond}) is satisfied.

 A special property of the operator $P=-\partial_x
a(x)\partial_x$ with step-like periodic coefficient $a(x)$ is that
the coefficients $(r\pm 1)$ in the dispersion relation
\begin{equation}\label{blochpar}
2\cos\theta(\lambda)=(\rho+1)\cos\{\lambda(x_2b_2
+(1-x_2)b_1)\}-(\rho-1)\cos\{\lambda(x_2b_2-(1-x_2)b_1)\}\end{equation}
are independent of the spectral parameter $\lambda.$ Formula
(\ref{blochpar}) implies that the band spectrum is periodic if the
profile of $a$ verifies (\ref{cond}).

The third part of the Theorem is proved in Section
\ref{s-convergence}.

The convergence of the resonances for a finitely periodic system
with $k$ cells  to the bands of real spectrum for the periodic
problem as $k\rightarrow\infty$ was discussed  by F. Barra and P.
Gaspard in \cite{BarraGaspard1999} in
 the case of Schr{\"o}dinger equation.

In our proof we use representations for the reflection and
transmission coefficients $r_k,$ $t_k$ for a finite slab of periodic
medium as in the recent paper of Molchanov and Vainberg
\cite{MolchanovVainberg2006}.
 The authors considered transition of truncated medium described by
the $1-$D Schr\"odinger operator to semi-infinite periodic
materials. By relating the reflection coefficients to the resolvent
of $P_k$ we show explicitly that the resonances correspond to the
poles of the analytic continuation of $r_k(\lambda)$ to $\cz_-.$
 Then we consider  the limit of the poles of $r_k(\lambda),$  as $k\rightarrow\infty.$

Note that for $\lambda\in\rz$ the reflection coefficient $r_{k+1}$
for $k+1$ cells medium is related to $r_k$ for $k$ cells medium via
$r_{k+1}=f_\lambda(r_k),$ where $f_{\lambda}$ is a linear-fractional
automorphism of the unit disk. By considering the fixed point of
$f_\lambda$ we get a new proof of the convergence of $r_k(\lambda)$
to $r(\lambda)$ when $\lambda$ belongs to the spectral gapes and
non-degenerate band edges for the operator $P$ (see Section
\ref{s-fixed_point}).


The structure of the paper is the following:\\
 In Section
\ref{s-molchanov} we recall some well-known facts concerning
spectral problem for weighted Sturm-Liouville operators (see
\cite{LevitanSargsyan}) and consider scattering by a finite slab of
a periodic medium. We follow \cite{MolchanovVainberg2006} with minor
changes due to the special form of  operator $P_k.$ We recall exact
formulas for the the reflection and transmission coefficients using
the iteration of the monodromy matrix. We recall also the result of
\cite{MolchanovVainberg2006} on a transition to semi-infinite
periodic material (limit $k\rightarrow\infty$). In Section
\ref{s-steplike} we give explicit expression for  the monodromy
matrix of $P_k.$ In Section \ref{s-banica} we recall the iterative
procedure used in \cite{Banica2003} for construction of the
resolvent $R_k$ and define resonances.   In Section \ref{s-relation}
the reflection coefficient $r_k$ is expressed using the iteration
formulas of \cite{Banica2003} and we show that the poles of
$R_k(\lambda)$ and the poles of $r_k(\lambda),$ $\lambda\in\cz_-,$
coincide.   In Section \ref{s-convergence} we prove the convergence
of ${\rm Res}\, (P_k)$ to the real axis. In Section
\ref{s-fixed_point} we discuss the convergence of $r_k(\lambda),$
$k\rightarrow\infty,$ for $\lambda\in\rz$ by considering  the limit
of a sequence of linear-fractional automorphisms on the unit disk.
In Appendix \ref{s-numerics} we present numerical examples.\\ \\
\textsc{Acknowledgements.}
  The author would like
to thank Maciej Zworski for suggesting to look at the problem
considered in the present paper and for helpful discussions.

\section{General methods for truncated periodic operators}\label{s-molchanov}
 In this section we following \cite{MolchanovVainberg2006} consider the scattering theory for operator $P_k$ combining the Floquet-Bloch theory and scattering theory for $1-$D weighted operators.

\subsection{The monodromy matrix and Bloch
quasi-momentum}\label{ss-Lyaponov} We recall first some well-known
facts
 concerning the spectral problem of Sturm-Liouville operators on the line (see \cite{LevitanSargsyan}). We consider equation
\begin{equation}\label{mainequation}
P\psi=-\partial_x a(x)\partial_x \psi (x)=\lambda^2\psi
\end{equation} on
$\{\psi\in H_{\rm loc}^1(\rz),\,\,a\partial_x\psi\in H_{\rm
loc}^1(\rz)\}$ with  a strictly positive $a (x)$ as in the
Introduction, formula (\ref{pk}) or periodic as $a_0$ in (\ref{a0}).

 Let $\psi_{1,2}$ be solutions of
(\ref{mainequation}) with initial data
\begin{equation}\label{initialdata}
\psi_1(\lambda,0)=1,\,\,(a\partial_x\psi_1)(\lambda,0)=0;\,\,\psi_2(\lambda,0)=0,\,\,(a\partial_x\psi_2)(\lambda,0)=1.
\end{equation}

We define the transfer matrix (propagator) $M_\lambda(0,x)$ for
operator $P$
\begin{equation}\label{pruffer}
M_\lambda(0,x)=\left(
                 \begin{array}{cc}
                   \psi_1 (\lambda, x) & \lambda\psi_2 (\lambda, x) \\
                   \frac{(a\partial_x\psi_1)(\lambda,x)}{\lambda} & (a\partial_x\psi_2)(\lambda ,x) \\
                 \end{array}
               \right).
               \end{equation}
From (\ref{initialdata}) it follows that $M_\lambda(0,x)$ is the
identity matrix. For any solution $\psi$ of (\ref{mainequation})
matrix $M_\lambda(0,x)$ maps the Cauchy data of $\psi$ at $x=0$ into
the Cauchy data of $\psi$ at point $x:$ $$M_\lambda(0,x):\,\,\left(
                       \begin{array}{c}
                         \psi(0) \\
                         \frac{(a\partial_x\psi)(\lambda,0)}{\lambda} \\
                       \end{array}
                     \right)\,\,\mapsto\,\,\left(
                       \begin{array}{c}
                         \psi(x) \\
                         \frac{(a\partial_x\psi)(\lambda,x)}{\lambda} \\
                       \end{array}
                     \right).$$
As the generalized Wronskian associated with $\psi_1,\psi_2$
$$W[\psi_1,\psi_2]=\psi_1 a\partial_x\psi_2-\psi_2
a\partial_x\psi_1$$ is constant, we have
$$\det M_\lambda(0,x)=W[\psi_1,\psi_2](1)=W[\psi_1,\psi_2](0)=1.$$
Equation (\ref{mainequation}) with $\Im\lambda >0$ has exactly one
solution $\psi^+$ in $L^2(\rz_+)$ normalized by the condition
$\psi^+(\lambda, 0)=1,$ and it has exactly one solution $\psi^-\in
L^2(\rz_-)$ normalized by the same condition. Here $\rz_\pm$ are
semiaxes $x\gtrless 0.$ Any solution of (\ref{mainequation}) can be
represented as linear combinations of $\psi_1$ and $\psi_2$ and from
the normalization of $\psi^\pm$ it follows, that there exist
functions $m^\pm=m^\pm (\lambda) $ such that
$$\psi^\pm=\psi_1 +m^\pm(\lambda)\psi_2,\,\,\Im\lambda >0.$$
Functions $m^\pm=m^\pm (\lambda)$ are called Weyl's functions and we
have
$$\psi_1 +m^+(\lambda)\psi_2\in L^2(\rz_+),\,\,\psi_1
+m^-(\lambda)\psi_2\in L^2(\rz_-),\,\,\Im\lambda >0.$$

Let $a(x)$ be periodic: $a(x+l)=a(x).$ Consider propagator through
one period (monodromy matrix):
$$M_\lambda=M_\lambda(0,l)=\left(
                             \begin{array}{cc}
                               \alpha & \beta \\
                               \gamma & \delta \\
                             \end{array}
                           \right)(\lambda)=\left(
                 \begin{array}{cc}
                   \psi_1 (\lambda, l) & \lambda\psi_2 (\lambda, l) \\
                   \frac{(a\partial_x\psi_1)(\lambda,l)}{\lambda} & (a\partial_x\psi_2)(\lambda ,l) \\
                 \end{array}
               \right).$$
Denote $F(\lambda)=\frac12 \tr M_\lambda =(\alpha +\delta)/2=(\psi_1
(\lambda,1)+a\partial\psi_2(\lambda,1)$ the Lyapunov function.
                           Both $M_\lambda$ and $F(\lambda)$ are entire function of
                           $\lambda$ and $\det M_\lambda =1.$
                           The eigenvalues $\mu^\pm(\lambda)$ of
                           $M_\lambda$ are the roots of the
                           characteristic equation \begin{equation}\label{chareq}
                           \mu^2 -2\mu F(\lambda) +1=0.\end{equation}
If $\Im\lambda >0$ then one can select roots $\mu^\pm (\lambda)$ of
(\ref{chareq}) in such a way that $\mu^\pm(\lambda)=e^{\pm
il\theta(\lambda)},$ where $\theta(\lambda)$ is analytic and
\begin{equation}\label{important}
\Im\theta(\lambda) >0\,\,\mbox{when}\,\,\Im\lambda >0,
\end{equation}
i.e., \begin{equation}\label{important2} |\mu^+(\lambda)|
<1,\,\,|\mu^-(\lambda)| >1,\,\,\Im\lambda >0.
\end{equation}
The roots $\mu^\pm (\lambda)$ for real $\lambda\geq 0$ are defined
by continuity in the upper half plane:
$$\mu^\pm (\lambda)=\mu^\pm (\lambda+i0),\,\,\lambda\in
[0,\infty).$$

Since the trace of $M_\lambda$ is equal to the sum of the
eigenvalues $e^{\pm il\theta(\lambda)}$,
   \begin{equation}\label{(12)}
   \cos l\theta =F(\lambda)=\frac12(\psi_1+a\psi_2')(\lambda,l)=\frac12(\alpha +\delta).
   \end{equation}

The spectrum of $P$ belongs to the positive part of the energy axis
$E=\lambda^2$ and  has band structure.

For real $\lambda \geq 0,$ the inequality $|F(\lambda)| \leq 1$
defines the spectral bands (zones)
$$b_n=[\lambda_{2n-1},\lambda_{2n}],\,\,n=1,2,\ldots,$$ on the
frequency axis $\lambda=\sqrt{E}$. The bands are defined by the
condition $|F(\lambda)| \leq 1$ and $F(\lambda)=\pm 1$ at any band
edge $\lambda=\lambda_j.$

The function $\theta(\lambda)$ is real valued when $\lambda$ belongs
to a band. The roots $\mu^\pm(\lambda)$ are complex adjoint there,
and $|\mu^\pm| =1.$

The spectrum of $P$ (on $L^2(\rz)$) on the frequency axis is
$\bigcup_{n=1}^\infty b_n.$

The complimentary open set, given by $|F(\lambda)| >1,$ corresponds
to spectral gaps, $\bigcup_{n=1}^\infty g_n.$ On gaps, the function
$i\theta(\lambda)$ is real valued,
 the roots
$\mu^\pm(\lambda)$ are real and (\ref{important2}) holds.

A point $\lambda_j$ which belongs to the boundary of a band and the
boundary of a gap is called a non-degenerate  band edge. If it
belongs to the boundary of two different bands, it is called a
degenerate band edge.

As in \cite{MolchanovVainberg2006} we get that if $\lambda
=\lambda_0$ is a non-degenerate band edge, then $F'(\lambda)\neq 0.$
If $\lambda=\lambda_0$ is a degenerate edge, then $F'(\lambda)=0,$
$F''(\lambda)\neq 0.$ Both eigenvalues of the monodromy matrix
$M_\lambda$ at any band edge are equal to $1$ or both are equal to
$-1.$

We normalize the eigenvectors $h^\pm(\lambda)$ of $M_\lambda$ by
choosing the first coordinate of $h^\pm (\lambda)$ to be equal to
one:
$$h^\pm (\lambda)=\left(
                    \begin{array}{c}
                      1 \\
                      m^\pm(\lambda) \\
                    \end{array}
                  \right).$$ The second coordinates of the vectors
                  $h^\pm (\lambda)$ coincide with the Weyl's
                  functions defined above.
                  In fact, if $\psi^\pm$ are solutions of the
                  equation $P\psi=\lambda^2\psi$ with the initial
                  Cauchy data given by the eigenvector $h^\pm,$ then
                  \begin{equation}\label{(13)}
                  \psi^\pm (\lambda, x+l)=e^{\pm il\theta (\lambda)}
                  \psi^\pm(\lambda,x) \end{equation}
                  and (\ref{important2}) implies that $\psi^\pm\in
                  L^2(\rz_\pm)$ when $\Im\lambda >0.$ From here it
                  follows that $\psi^\pm$ coincide with Weyl's
                  solution introduced for general Hamiltonians $P,$
                  and that the second coordinates of the vectors
                  $h^\pm(\lambda)$ are Weyl's functions.

Since $$\left(
          \begin{array}{cc}
            \alpha-e^{\pm il \theta(\lambda)} & \beta \\
            \gamma & \delta-e^{\pm il \theta(\lambda)} \\
          \end{array}
        \right)\left(
                 \begin{array}{c}
                   1 \\
                   m^\pm \\
                 \end{array}
               \right)=0,$$
the following two representations are valid for Weyl's functions:

\begin{equation}\label{representations}
m^\pm(\lambda)=\frac{e^{\pm il\theta(\lambda)}-\alpha
(\lambda)}{\beta(\lambda)}=\frac{\gamma(\lambda)}{e^{\pm
il\theta(\lambda)} -\delta(\lambda)}.
\end{equation}

\subsection{Reflection coefficient for the truncated periodic operator}\label{ss-transition} We consider operator $P_k$ with the
truncated periodic coefficient $a_k:$
$$P_k\psi =-\partial_x a_k(x)\partial_x\psi
(x),\,\,a_k(x)=a(x),\,\,\mbox{for}\,\,x\in
[0,kl];\,\,a_k(x)=\frac{1}{b_1^2},\,\,\mbox{for}\,\,x\not\in [0,kl]
,$$ which appears when the propagation of waves through a finite
slab of a periodic medium is studied. We shall also consider the
limiting case $k=\infty:$
$$P_\infty\psi =-\partial_x a_\infty(x)\partial_x\psi (x),$$ which
corresponds to the case of such a long slab that it can be
considered as half infinite.

We shall denote the reflection and transmission coefficients for the
operator $P_k$ (with compactly supported coefficient $a_k$) by $r_k$
and $t_k,$ respectively:
$$P_k\psi =\lambda^2\psi,\,\,\psi=e^{i\lambda b_1x} +r_k e^{-i\lambda
b_1 x},\,\,x <0;\,\,\psi=t_ke^{i\lambda b_1x} ,\,\,x >kl.$$

In the case of the operator $P_\infty,$ the solution $\psi$ of the
scattering problem is defined as the solution of the equation
$P_\infty\psi =\lambda^2\psi,$ such that
\begin{equation}\label{(19)}\psi=e^{i\lambda b_1x} +r e^{-i\lambda b_1 x},\,\,x
<0;\,\,\psi= c \psi^+(\lambda, x) ,\,\,x
>0,\end{equation} with some $r=r(\lambda),$ $c=c(\lambda).$
We have the following version of  Theorem 3 of S. Molchanov, B.
Vainberg in \cite{MolchanovVainberg2006}:
\begin{theorem}\label{tMV}
{\bf 1)} The transfer matrix over $k$ periods $M^k_\lambda=
T_\lambda(0,lk)$ has the form
\begin{equation}\label{(20)} M_\lambda^k=\left(
                                           \begin{array}{cc}
                                             \alpha_k & \beta_k \\

                                             \gamma_k & \delta_k \\
                                           \end{array}
                                         \right)=\frac{\sin
                                         kl\theta(\lambda)}{\sin
                                         l\theta(\lambda)} M_\lambda
                                         -\frac{\sin (k-1)l\theta
                                         (\lambda)}{\sin
                                         l\theta(\lambda)} I,
\end{equation} where $\theta=\theta(\lambda)$ is the Bloch function.
The elements of $M_\lambda^k$ satisfy the relations
\begin{align}
&\alpha_k-\delta_k=\frac{\sin kl\theta(\lambda)}{\sin
l\theta(\lambda)}(\alpha -\delta),\,\,\,\,\beta_N=\frac{\sin
kl\theta(\lambda)}{\sin l\theta(\lambda)}\beta,\nonumber\\
&\gamma_k=\frac{\sin kl\theta(\lambda)}{\sin
l\theta(\lambda)}\gamma,\,\,\,\, \alpha_k+\delta_k=2\cos
kl\theta(\lambda).\label{(21)}
\end{align}
{\bf 2)} The reflection coefficients have the forms
\begin{equation}\label{(23)}
r_k(\lambda)=-\frac{(\alpha-\delta)+i(b_1\gamma +\frac{\beta}{b_1}
)}{2\sin l\theta (\lambda)\frac{\cos kl\theta(\lambda)}{\sin
kl\theta(\lambda)} +i(b_1\gamma -\frac{\beta}{b_1})},
\end{equation}
\begin{equation}\label{(24)}
r(\lambda)=\frac{\frac{\beta}{b_1} +b_1\gamma
-i(\alpha-\delta)}{2\sin l\theta (\lambda) -(b_1\gamma
-\frac{\beta}{b_1})}.
\end{equation}
{\bf 3)} The transmission probability have the form
\begin{equation}\label{(25)}
|t_k(\lambda)|^2=\frac{4}{\frac{\sin^2
lk\theta}{\sin^2l\theta}\left((\alpha -\delta )^2 +(b_1\gamma
+\frac{\beta}{b_1})^2\right)+4}=\frac{1}{\frac{\sin^2
lk\theta}{\sin^2l\theta}\frac{|r_1|^2}{|t_1|^2}+1}.
\end{equation}
{\bf 4)}The reflection coefficients $r_k(\lambda)$ and $r(\lambda)$
are analytic in the upper half plane $\cz_+=\{\lambda:\,\,\Im\lambda
>0\}$ and continuous in $\overline{\cz}_+.$
For any $\lambda\in \overline{\cz}_+\setminus \cup_{n=1}^\infty b_n$
we have $r_k(\lambda)\rightarrow r(\lambda).$ When
$\lambda\in\cup_{n=1}^\infty b_n,$  $r_k(\lambda)$ converges to
$r(\lambda)$ in the weak sense:
$$
\int_{-\infty}^\infty r_k(\lambda)\varphi (\lambda) d\lambda
\rightarrow \int_{-\infty}^\infty r(\lambda)\varphi (\lambda)
d\lambda\,\,\,\,\mbox{as}\,\, k\rightarrow \infty.$$ for any test
function $\varphi\in D.$
\end{theorem}
{\bf Proof:}  We reproduce here the proof of
\cite{MolchanovVainberg2006} for the sake of completeness with only
minor changes due to the ``weight'' in the  definitions of $P_k$ and
$P_\infty.$  Formula (\ref{(20)}) follows by induction from relation
(\ref{chareq}):
$$M_\lambda^2 -2\cos l\theta (\lambda) M_\lambda +I =0.$$
The first three relations of (\ref{(21)}) are immediate consequences
of (\ref{(20)}). In order to get the fourth one we note that the
eigenvalues of $M_\lambda^k$ are $\mu^\pm(\lambda)=e^{\pm ikl\theta
(\lambda)}.$ Thus,
\begin{equation}\label{(26)}
\alpha_k +\delta_k=\tr M_\lambda^k= e^{ikl\theta(\lambda)}
+e^{-ikl\theta (\lambda)} =2\cos kl\theta (\lambda).
\end{equation}

Next we prove (\ref{(23)}). The relation between Cauchy data for the
left-to-right scattering solution at $x=0$ and $x=kl$ are given by
$$\left(
    \begin{array}{cc}
      \alpha_k & \beta_k \\
      \gamma_k & \delta_k \\
    \end{array}
  \right)\left(
           \begin{array}{c}
             1+r_k \\
             \frac{i}{b_1}(1-r_k) \\
           \end{array}
         \right)=\left(
                   \begin{array}{c}
                     t_k \\
                     \frac{it_k}{b_1} \\
                   \end{array}
                 \right)\,\,\Leftrightarrow\,\,\left\{\begin{array}{c}
                                                 \alpha_k(1+r_k)
                 +\beta_k(\frac{i}{b_1} (1-r_k)) =t_k \\
                                                 \gamma_k (1+r_k)
                                                 +\delta_k
                                                 (\frac{i}{b_1}(1-r_k))=\frac{it_k}{b_1}
                                               \end{array}\right.
$$ By dividing the second equation by the first one we arrive at
$$\frac{\gamma_k (1+r_k) +\frac{\delta_k i}{b_1}(1-r_k)}{\alpha_k
(1+r_k)+\frac{\beta_k i}{b_1}(1-r_k)}=\frac{i}{b_1}.$$ Solving for
$r_k$ we obtain \begin{equation}\label{rk} r_k=\frac{\delta_k
-\alpha_k -i(b_1\gamma_k +\frac{\beta_k}{b_1})}{\delta_k +\alpha_k
+i(b_1\gamma_k -\frac{\beta_k}{b_1})}.\end{equation} Using
(\ref{(21)}) we get
$$r_k=-\frac{\frac{\sin kl\theta}{\sin l\theta} (\alpha
-\delta)+i\left( b_1\frac{\sin kl\theta}{\sin l\theta}\gamma
+\frac{1}{b_1}\frac{\sin kl\theta}{\sin l\theta}\beta\right)}{2\cos
kl\theta +i\left( b_1\frac{\sin kl\theta}{\sin l\theta}\gamma
-\frac{1}{b_1}\frac{\sin kl\theta}{\sin l\theta}\beta\right)}.$$
This justifies (\ref{(23)}).

In order to get (\ref{(24)}) we note that (\ref{(19)}) implies that
$$\left(
           \begin{array}{c}
             1+r \\
             \frac{i}{b_1}(1-r)\\
           \end{array}
         \right)=c\left(\begin{array}{c}
                    1 \\
                    m^+
                  \end{array}\right)
                  \,\,\Leftrightarrow\,\,\frac{1+r}{\frac{i}{b_1}(1-r)}=\frac{1}{m^+}$$
                  and therefore, $$r=\frac{\frac{i}{b_1}-m^+}{\frac{i}{b_1}+m^+}.$$
From here and (\ref{representations}) it follows that
$$r=\frac{\frac{i}{b_1}-\frac{e^{il\theta} -\alpha}{\beta}}{\frac{i}{b_1}+\frac{e^{il\theta} -\alpha}{\beta}}=\frac{\alpha +\frac{i}{b_1}\beta -e^{il\theta}}{e^{il\theta}-(\alpha -\frac{i}{b_1}\beta)}.$$
and
$$r=\frac{\frac{i}{b_1}-\frac{\gamma}{e^{il\theta}-\delta}}{\frac{i}{b_1}+\frac{\gamma}{e^{il\theta}-\delta}}=\frac{\frac{i}{b_1}(e^{il\theta}-\delta)-\gamma}{\frac{i}{b_1}(e^{il\theta}-\delta)+\gamma}=
\frac{e^{il\theta}-(\delta-ib_1\gamma)}{e^{il\theta}-(\delta+ib_1\gamma)}.$$
Hence,
$$r\left(2 e^{il\theta}-(\alpha -\frac{i}{b_1}\beta) -(\delta +i
b_1\gamma)\right)=\alpha +\frac{i}{b_1}\beta -(\delta -ib_1\gamma)$$
and \begin{align*}r&=\frac{(\alpha +\frac{i}{b_1}\beta) -(\delta
-ib_1\gamma)}{2 e^{il\theta}-(\alpha -\frac{i}{b_1}\beta) -(\delta
+i b_1\gamma)}=\frac{(\alpha-\delta)
+i(\frac{\beta}{b_1}+b_1\gamma)}{2
e^{il\theta}-(\alpha+\delta)+i(\frac{\beta}{b_1}-b_1\gamma)}=\\
&=\frac{(\alpha-\delta) +i(\frac{\beta}{b_1}+b_1\gamma)}{ 2i\sin
l\theta+i(\frac{\beta}{b_1}-b_1\gamma)},
\end{align*}
where the last equality is a consequence of (\ref{(12)}) and it
implies (\ref{(24)}).

We prove the third statement of the theorem. From (\ref{rk}) it
follows \begin{equation}\label{rk2}|r_k|^2=\frac{(\delta_k
-\alpha_k)^2 + (b_1\gamma_k +\frac{\beta_k}{b_1})^2}{(\delta_k
+\alpha_k)^2 + (b_1\gamma_k -\frac{\beta_k}{b_1})^2}.\end{equation}
Using that $|r_k|^2+|t_k|^2=1$ we get
\begin{equation}\label{tk2}|t_k|^2=1-|r_k|^2=\frac{4\delta_k\alpha_k
-4\gamma_k\beta_k}{(\delta_k +\alpha_k)^2 + (b_1\gamma_k
-\frac{\beta_k}{b_1})^2}.\end{equation} We use $\det
M_\lambda^k=\alpha_k\delta_k-\beta_k\gamma_k=1$ and get
$$|t_k|^2=\frac{4}{4+(\delta_k-\alpha_k)^2 +(b_1\gamma_k
+\frac{\beta_k}{b_1})^2}=\frac{4}{\frac{\sin^2
lk\theta}{\sin^2l\theta}\left((\alpha -\delta )^2 +(b_1\gamma
+\frac{\beta}{b_1})^2\right)+4}.$$

From formulas (\ref{rk2}) and (\ref{tk2}) we get
$$\frac{|r_k|^2}{|t_k|^2}=\frac14\left((\delta_k -\alpha_k)^2 + (b_1\gamma_k
+\frac{\beta_k}{b_1})^2\right)$$ and hence, putting $k=1,$
$$|t_k|^2=\frac{1}{\frac{\sin^2
lk\theta}{\sin^2l\theta}\frac{|r_1|^2}{|t_1|^2}+1}.$$

The analyticity of $r_k(\lambda)$ and $r(\lambda)$ in $\cz_+$ and
their continuity in $\overline{\cz}_+$ follow from the explicit
formulas (\ref{(23)}), (\ref{(24)}). For
$\lambda\in\cz_+\bigcup\cup_{n=1}^\infty g_n,$ we have
$\Im\theta(\lambda) >0$ and $\theta(\lambda)$ is pure imaginary on
the gaps $g_n.$ Furthermore, if $\Im\theta (\lambda) >0,$ then
$$\frac{\cos k\theta(\lambda)}{\sin k\theta (\lambda)}\rightarrow -i
\,\,\mbox{as}\,\, k\rightarrow \infty $$ and this justifies the
convergence of $r_k(\lambda)$ to $r$ when
$\lambda\in\cz_+\bigcup\cup_{n=1}^\infty g_n.$

The weak convergence for $\lambda\in\cup_{n=1}^\infty b_n$  is a
consequence the convergence in the complex half plane.

The proof of Theorem \ref{tMV} is complete.
 \hfill\qed

Note also the following relations:
$$\frac{|r_k|^2}{|t_k|^2}=\frac{\sin^2
lk\theta}{\sin^2l\theta}\cdot\frac{|r_1|^2}{|t_1|^2},$$
\begin{align*}|t_k|^2&=\frac{4}{4\cos^2 kl\theta +\left( b_1\frac{\sin
kl\theta}{\sin l\theta}\gamma -\frac{1}{b_1}\frac{\sin
kl\theta}{\sin l\theta}\beta\right)^2}=\frac{4\frac{\sin^2
l\theta}{\sin^2 kl\theta}}{4\sin^2 l\theta\frac{\cos^2
kl\theta}{\sin^2 kl\theta} +\left( b_1\gamma
-\frac{\beta}{b_1}\right)^2}.\end{align*} The last formula follows
from (\ref{tk2}) by using   (\ref{(26)}).

Formula (\ref{(25)}) implies that the perfect transmission
($|t_k|^2=1$) occurs whenever $|r_1|^2=0$ ($ |t_1|^2=1$) or if
\begin{equation}\label{eek}
\frac{\sin^2 lk\theta}{\sin^2l\theta}=0.\end{equation} For
$\theta\in [0,\pi/l]$ equation (\ref{eek}) is satisfied when $\theta
lk=m\pi$ for $m=1,2,\ldots,k-1.$

 Therefore,
in the general case ($|r_1|^2\neq 0$), the transmission probability
has $k-1$ peaks with $|t_k|^2=1$ in each allowed energy band as
$\theta$ increases by $\pi/l.$ Since the peaks in the transmission
probability (or in general in the cross section) are associated with
resonances, we expect to find  $k-1$ resonances near each allowed
energy band.

On the gaps, the function $i\theta(\lambda)$ is real valued. Then
the transmission probability is given by
\begin{equation}\label{(25gaps)}
|t_k(\lambda)|^2=\frac{1}{\frac{\sinh^2
lki\theta}{\sinh^2li\theta}\frac{|r_1|^2}{|t_1|^2}+1}.
\end{equation}
As  $\sinh^2 lki\theta\neq 0$ for $\theta\neq 0,$ then in the
forbidden zone $|t_k(\lambda)|^2\neq 1$ unless $|t_1|^2=1,\,\,
|r_1|^2=0.$ Hence there are no resonances below the gaps.

On the gaps,  the reflection coefficient for the half-periodic
system $r(\lambda)$ satisfy
$$|r(\lambda)|^2=\frac{\left(\frac{\beta}{b_1}
+b_1\gamma\right)^2+(\alpha -\delta)^2}{e^{2i\theta}+
e^{-2i\theta}-2 +\left( b_1\gamma
-\frac{\beta}{b_1}\right)^2}=\frac{\left(\frac{\beta}{b_1}
+b_1\gamma\right)^2+(\alpha -\delta)^2}{(\alpha +\delta)^2-4 +\left(
b_1\gamma -\frac{\beta}{b_1}\right)^2}=1.$$



In \cite{MolchanovVainberg2006}, Theorem 5, was shown that
\begin{lemma}
If $\lambda_0$ is a degenerate band edge, i.e. $F(\lambda_0)=\pm
1,\,\,F'(\lambda_0)=0,$ then the reflection coefficient $r_k$ is
zero, $r_k(\lambda_0)=0.$
\end{lemma}

 The proof uses the fact that at any degenerate
band edge $\lambda=\lambda_0$ the monodromy matrix
$M_{\lambda_0}=\pm I$ and $F''(\lambda_0)\neq 0.$ This allows to
pass to the limit in (\ref{(23)}) as $\lambda\rightarrow\lambda_0.$
The numerator in the \rhs{} of (\ref{(23)}) vanishes as
$\lambda\rightarrow\lambda_0.$ The denominator converges to $\pm
2/k,$ since $\theta(\lambda_0)=n\pi.$

Thus the medium is transparent for the plane wave with the frequency
$\lambda_0$ and we expect to find a resonance $\lambda\in\cz_-$
below the degenerated band edge $\lambda_0.$

\section{The monodromy matrix for $P_k$}\label{s-steplike}

In this section we give expressions for the elements of the
monodromy matrix $M_\lambda.$

Let $a(x)$ be $1-$periodic function equal to $a_0$ for $x\in [0,1)$
as in (\ref{pk}), (\ref{a0}):
\begin{equation*}
a_0 (x)=\left\{\begin{array}{lr}
                   b_2^{-2} & \mbox{for}\,\,x\in [0, x_2) \\
                   b_1^{-2} & \mbox{for}\,\, x\in [x_2,1)
                 \end{array}\right.\end{equation*} and

\begin{equation*}a (x)=\left\{\begin{array}{lc}
                   a_0(x-j) & \mbox{for}\,\,x\in [j, j+1),\,\,0\leq j\leq k-1 \\
                   b_1^{-2} & \mbox{elsewhere}
                 \end{array}\right.\end{equation*}
                 for $k\geq 2,$  where $k$ is the number of identical
                 cells. The period $l=1.$
The normalized solutions $\psi_1(\lambda,x),$ $\psi_2(\lambda,x),$
$$-\partial_x
a(x)\partial_x\psi_i=\lambda^2\psi_i,\,\,\psi_1(\lambda,0)=(b_2^{-2}\partial_x\psi_2)(\lambda,0)=1,\,\,\psi_2(\lambda,0)=(b_2^{-2}\partial_x\psi_1)(\lambda,0)=0,$$
satisfy (see \cite{Banica2003})
\begin{align*}
&\psi_1(\lambda,x)=\left\{\begin{array}{cc}
                     \frac12 e^{i\lambda b_2 x} +\frac12 e^{-i\lambda b_2 x} & \mbox{for}\,\,x\in [0,x_2) \\
                     A e^{i\lambda b_1 x} +Be^{-i\lambda b_1 x} &
                     \mbox{for}\,\, x\in [x_2,1)
                   \end{array}\right.\\
&\psi_2(\lambda,x)=\left\{\begin{array}{cc}
                     -\frac{ib_2}{2\lambda} e^{i\lambda b_2 x} +\frac{ib_2}{2\lambda} e^{-i\lambda b_2 x} & \mbox{for}\,\,x\in [0,x_2) \\
                     C e^{i\lambda b_1 x} +De^{-i\lambda b_1 x} &
                     \mbox{for}\,\, x\in [x_2,1)
                   \end{array}\right.
\end{align*}
with $A,B,C,D$ chosen such that $\psi_j$ and $a(x)\partial_x\psi_j$
are continuous at $x_2:$
\begin{align*}
& A=\frac{1}{4b_2}\left[ (b_2+b_1)e^{i\lambda x_2 (b_2-b_1)}
+(b_2-b_1)e^{-i\lambda x_2 (b_2 +b_1)}\right],\\
& B=\frac{1}{4b_2}\left[ (b_2+b_1)e^{-i\lambda x_2 (b_2-b_1)}
+(b_2-b_1)e^{i\lambda x_2 (b_2 +b_1)}\right],\\
&C=\frac{i}{4\lambda}\left[ -(b_2+b_1)e^{i\lambda x_2
(b_2-b_1)} +(b_2-b_1)e^{-i\lambda x_2 (b_2 +b_1)}\right],\\
&D=\frac{i}{4\lambda}\left[ (b_2+b_1)e^{-i\lambda x_2 (b_2-b_1)}
-(b_2-b_1)e^{i\lambda x_2 (b_2 +b_1)}\right].
\end{align*}

 We get the monodromy matrix
$$M_\lambda=\left(
              \begin{array}{cc}
                \alpha & \beta \\
                \gamma & \delta \\
              \end{array}
            \right)=\left(
              \begin{array}{cc}
                \psi_1(\lambda,1) & \lambda\psi_2(\lambda,1) \\
                \frac{1}{\lambda}(\frac{1}{b_1^2}\partial_x\psi_1)(\lambda,1) & (\frac{1}{b_1^2}\partial_x\psi_2)(\lambda,1) \\
              \end{array}
            \right)
$$ with
\begin{align*}
&\alpha=\psi_1(\lambda,1)=\frac{b_2+b_1}{2b_2}\cos\lambda
[b_1(1-x_2) +x_2b_2] +\frac{b_2-b_1}{2b_2}\cos\lambda [b_1(1-x_2)
-x_2b_2],\\
&\beta=\lambda\psi_2(\lambda,1)=\frac{b_2+b_1}{2}\sin\lambda
[b_1(1-x_2)+x_2b_2]-\frac{b_2-b_1}{2}\sin\lambda
[b_1(1-x_2)-x_2b_2],\\
&\gamma =\frac{1}{\lambda}(\frac{1}{b_1^2}\partial_x
\psi_1)(\lambda, 1)=-\frac{b_2+b_1}{2b_1b_2}\sin\lambda
[b_1(1-x_2)+x_2b_2]
-\frac{b_2-b_1}{2b_1b_2}\sin\lambda [b_1(1-x_2)-x_2b_2],\\
&\delta=(\frac{1}{b_1^2}\partial_x\psi_2)(\lambda,1)=\frac{b_2+b_1}{2b_1}\cos\lambda
[b_1(1-x_2)+x_2b_2] - \frac{b_2-b_1}{2b_1}\cos\lambda
[b_1(1-x_2)-x_2b_2].
\end{align*}
Then $$\tr(M_\lambda)=\alpha+\delta=(\rho+1)\cos\lambda
[b_1(1-x_2)+x_2b_2]-(\rho-1)\cos\lambda [b_1(1-x_2)-x_2b_2],$$ where
$\rho=\frac{b_2^2+b_1^2}{2b_2b_1}.$ The Bloch quasi-momentum
$\theta=\theta(\lambda)$ satisfy
$$2\cos\theta=\tr M_\lambda.$$





The first formula in (\ref{(25)}) implies that the one cell medium
is perfectly transparent:
$|t_1(\lambda_0)|^2=1,\,\,|r_1(\lambda_0)|^2=0,$ if
$\alpha-\delta=0,$ $b_1\gamma +\frac{\beta}{b_1}=0 .$ We get
$\lambda=\lambda_0=\pi m/ x_2b_2.$


Note that there is a resonance $\lambda\in{\rm Res}\,(P_1)$ for the
one cell operator such that $\Re\lambda=\lambda_0$ (see equation
(\ref{n3}) with $b_3=b_1$ and $x_1=0$).

 If $\lambda_0 =\frac{\pi
m}{x_2 b_2}$  then we have
\begin{align*}2F(\lambda_0)=2\cos\theta (\lambda_0)&= \pm2\cos
\left(\frac{(1-x_2)b_1}{x_2b_2}\pi m \right),\,\,\mbox{if $m$ is
even (odd)}.
\end{align*}

Thus in general situation, $\frac{(1-x_2)b_1}{x_2b_2}\not\in\qz,$
non-rational, $\lambda_0$ is an interior point of a spectral band.

If $x_2b_2=(1-x_2)b_1$ then  the Lyapunov function satisfies
\begin{equation}\label{Lyp}
2F(\lambda)=(\rho+1)\cos(2\lambda b_2 x_2) -(\rho-1)\end{equation}
and we have
 $F(\lambda_0)=1,$ $F'(\lambda_0)=0$ and $F''(\lambda_0)\neq 0.$
 Hence
$\lambda_0$ is degenerate band edge.

The non-degenerate band edge is then given by the equation
$$(\rho+1)\cos(2\lambda b_2 x_2)
-(\rho-1)=-2\,\,\Leftrightarrow\,\,
\lambda=\frac{1}{2b_2x_2}\left(\pm\arccos
\left(\frac{\rho-3}{\rho+1}\right) +2\pi n\right).$$

\section{Explicit construction of the resolvent and resonances}\label{s-banica}
In this section we define the resonances as the poles of the
analytic continuation of the resolvent $R(\lambda)$  to $\cz_-.$
\subsection{Representation of the resolvent}



In this section we revue some formulas used by Valeria Banica in
\cite{Banica2003}, where she considered  the local and global
dispersion and the Stricharts inequalities for certain
one-dimensional Schr\"odinger and wave equations with step-like
coefficients. The systems  are described by the one-dimensional
Schr\"odinger equation
\begin{equation}\label{schr-eq} \left\{\begin{array}{rl}
     (i\partial_t +\partial_x a(x)\partial_x)u(t,x) &=0\,\,\mbox{for}\,\,(t,x)\in(0,\infty)\times\rz, \\
    u(0,x) & =u_0(x)\in L^2(\rz) \\
  \end{array}\right.
\end{equation}
    or by the one-dimensional wave equation
 \begin{equation}\label{wave-eq} \left\{\begin{array}{rl}
     (\partial_t^2 -\partial_x a(x)\partial_x)v(t,x) &=0\,\,\mbox{for}\,\,(t,x)\in\rz\times\rz, \\
    v(0,x) & =u_0(x)\in L^2(\rz), \\
    \partial_tv(0,x)&=0
  \end{array}\right.
\end{equation}
for a positive step-like function $a(x)$ with a finite number of
discontinuities.

 Consider a partition of the real axis
$$-\infty =x_0 <x_1 <x_2 <\ldots <x_{n-1} <x_n=\infty$$ and a step
function  $$a(x) =b_i^{-2}\,\,\mbox{for}\,\,x\in (x_{i-1},x_i),$$
where $b_i$ are positive numbers.

The operator $P:=-\partial_x a(x)\partial_x$ defined from $\{ u\in
H^1(\rz),\,\,a\partial_x u\in H^1(\rz)\}$ to $L^2(\rz)$ is
self-adjoint. For $\Im\lambda > 0,$  we define $R(\lambda)
v=(-\partial_x a(x)\partial_x -\lambda^{2})^{-1} v$ the resolvent of
$P.$

  Our choice of the spectral parameter $\lambda$
is related to $\omega$ in \cite{Banica2003} by  $\lambda= i\omega.$
We use the expression for the resolvent obtained in
\cite{Banica2003}.

On each interval $(x_i,x_{i+1})$ $R(\lambda) v$ is a finite sum of
terms
\begin{equation}\label{resolvent}
R(\lambda)v(x)=-\sum_{\mbox{finite}} C
e^{i\lambda\beta(x)}\int_{I(x_i)}\frac{v(y)}{2i\lambda}\frac{e^{\pm
i\lambda b_i y}}{\det D_n(-i\lambda)} dy
-\int_{-\infty}^\infty\frac{v(y)}{2i\lambda} b_i e^{i\lambda b_i
|x-y|} dy,\end{equation} where $I(x_i)$ is either $(-\infty, x_i)$
or $(x_i,\infty)$ and
$$\det D_2(-i\lambda)=(b_1+b_2)e^{i\lambda x_1(b_2-b_1)}.$$
In \cite{Banica2003} Banica defines all $\det D_n$ by induction. Let
$$\det \widetilde{D_2(-i\lambda)}=(b_2-b_1)e^{-i\lambda x_1
(b_1+b_2)}.$$ We have the following induction relations for $n\geq
3$
\begin{align}&\det D_n=e^{i\lambda b_n x_{n-1}}\left[(b_{n-1} -b_n)
e^{i\lambda b_{n-1} x_{n-1}}\det\widetilde{D_{n-1}}-(b_{n-1} +b_n)
e^{-i\lambda b_{n-1}x_{n-1}}\det D_{n-1}\right]\label{indrel}\\
&\det \widetilde{D_n}=e^{-i\lambda b_n x_{n-1}}\left[(b_{n-1} -b_n)
e^{-i\lambda b_{n-1} x_{n-1}}\det D_{n-1}-(b_{n-1} +b_n) e^{i\lambda
b_{n-1}x_{n-1}}\det \widetilde{ D_{n-1}}\right] .\nonumber
\end{align}

We define for $n\geq m\geq 2$
$$Q_{m}(-i\lambda)=e^{2i\lambda b_mx_m}\frac{\det
\widetilde{D_m}}{\det D_m},\,\, d_{m-1}=\frac{b_{m-1}-b_m}{b_{m-1}
+b_m}.$$ Then we have for $n\geq 3$ $$\det
D_n(i\lambda)=(b_1+b_2)e^{i\lambda (b_2-b_1)x_1}\Pi_{j=2\ldots
n-1}(b_j+b_{j+1})e^{-i\lambda
(b_j-b_{j+1})x_j}(1-d_jQ_j(-i\lambda)).$$

We have induction formula on the $Q_m$'s
\begin{equation}\label{linfr}
Q_m(-i\lambda)=e^{2i\lambda b_m (x_m-x_{m-1})}\frac{-d_{m-1}
+Q_{m-1}(-i\lambda)}{1-d_{m-1}Q_{m-1}(-i\lambda)}.\end{equation}


Note that a  linear-fractional transform on the unit disc occurs in
(\ref{linfr}) for $\Im\lambda \geq 0.$

If $\Im\lambda \leq 0$ then $|e^{2i\lambda b_m (x_m-x_{m-1})}|\geq
1. $ We use that $|d_n| <1$ and for any $n$ we can find $\epsilon_n
>0$  such that for every complex $\lambda$ with
$$\Im\lambda \geq -\epsilon_n$$ the estimate
$$|Q_2(-i\lambda)|=|d_1 e^{2i\lambda b_2(x_2-x_1)}|<1$$ holds
and gives by induction $$|Q_m(-i\lambda)| <1,\,\,n\geq m\geq 2.$$
Hence $(\det D_n(-i\lambda))^{-1}$ is uniformly bounded and well
defined in this region, which contains the real axis. Therefore
$i\lambda R(\lambda)u_0(x)$ can be analytically continued. The
spectral theorem gives
\begin{lemma} \label{l-1}The solution of the Schr\"odinger equation
(\ref{schr-eq}) verifies $$u(t,x)=\int_{-\infty}^\infty
e^{it\lambda^2}\lambda R(\lambda)u_0(x)\frac{d\lambda}{\pi}.
$$

The solution of the wave equation (\ref{wave-eq}) verifies
$$v(t,x)=\int_{-\infty}^\infty
e^{it\lambda}i\lambda R(\lambda)u_0(x)\frac{d\lambda}{2\pi}.
$$
\end{lemma}


Due to formula (\ref{resolvent}), the resonance spectrum ${\rm
Res}\,(P)$ consist of    zeros of $\det D_n(-i\lambda)$ or
equivalently the zeros of $1-d_{n-1}Q_{n-1}(-i\lambda)$ as by
(\ref{indrel}), equation $\det D_n(-i\lambda)=0$ is equivalent to
$$ Q_{n-1}(-i\lambda)=\frac{1}{d_{n-1}}.$$ By considerations
before Lemma \ref{l-1}, for each $n\geq 2$ there is $\epsilon_n>0$
such that all resonances verify $\Im\lambda <-\epsilon_n.$
 For $n=2$ there are no zeros.

For $n=3$ the resonances are  solutions of the equation $\det
D_3(-i\lambda)=0,$
$\lambda=\lambda_1+i\lambda_2,\,\,\lambda_j\in\rz,$ with constant
imaginary part:
\begin{align}
&\lambda_1=\frac{\pi
m}{2b_2(x_2-x_1)},\,\,\lambda_2=-\frac{1}{2b_2(x_1-x_2)}\ln\left|\frac{(b_2-b_3)(b_2-b_1)}{(b_2+b_3)(b_1+b_2)}
\right|<0,\label{n3}
\end{align}
where $m=0,\pm 2,\pm 4\,\ldots,$ is even, if $(b_2-b_3)(b_2-b_1)>0$
or $m=\pm 1,\pm 3,\pm 5\,\ldots,$ is odd, if $(b_2-b_3)(b_2-b_1)<0.$

 For $n>3$
the resonance spectrum can be obtained numerically using the
induction relations (\ref{indrel}).

\subsection{Locally periodic media}\label{ss-period}

 Suppose that the profile of $a(x)$ consist of a finite number
identical elements, obtained  by juxtaposing of $k$ unit cells.
Outside the interval $[0,k)$ the coefficient is constant. We make
the following choice: suppose $n=2k+1$ odd,
$b_1=b_3=b_5=\ldots=b_n,$ $b_2=b_4=\ldots=b_{n-1},$ and put
$$x_1=0,\,\, 0<x_2
<1,\,\,x_3=1,\,d=\frac{b_2-b_1}{b_2+b_1},\,\,\lambda=\lambda_1+i\lambda_2,\,\,\lambda_1\in\rz,\,\,\lambda_2
<0.$$ We have $x_{2k+1}=k,\,\,x_{2k} =k- 1+x_2.$  Let
\begin{equation}\label{a00}
a_0 (x)=\left\{\begin{array}{lr}
                   b_2^{-2} & \mbox{for}\,\,x\in [0, x_2) \\
                   b_1^{-2} & \mbox{for}\,\, x\in [x_2,1)
                 \end{array}\right.\end{equation} and

\begin{equation}\label{a}a (x)=\left\{\begin{array}{lc}
                   a_0(x-j) & \mbox{for}\,\,x\in [j, j+1),\,\,0\leq j\leq k-1 \\
                   b_1^{-2} & \mbox{elsewhere}
                 \end{array}\right.\end{equation}
                 for $k\geq 2.$  Here $k$ is the number of identical cells.

                The function $a(x)$ is called {\em locally periodic} or
                 {\em finite periodic} on the interval $[0,k)$ with $k$ cells.

 Then, with $n=2k+1,$
\begin{align*} &\det
D_n(-i\lambda)= (b_1+b_2)^{n-1}e^{-i\lambda(b_2-b_1)kx_2}
\Pi_{i=2,\ldots,n-1}\left(
1-(-1)^idQ_i(-i\lambda)\right),\end{align*} where
\begin{equation}\label{Qper}
Q_{2k} (-i\lambda)=e^{2i\lambda b_2 (k-1+x_2))}\frac{\det
\tilde{D}_{2k}}{\det D_{2k}},\,\,Q_{2k+1} (-i\lambda)=e^{2i\lambda
b_1 k}\frac{\det \tilde{D}_{2k+1}}{\det D_{2k+1}}
\end{equation}

 For $k=1,2,3,\ldots$ the
resonances are the solutions of the equation
\begin{equation}\label{detres}
\det D_{2k+1}(-i\lambda)=0
\end{equation} or equivalently
\begin{equation}\label{maineq}
Q_{2k}(-i\lambda)=\frac{1}{d},
\end{equation}
 where
\begin{align}
&Q_{2k}=e^{2i\lambda b_2x_2}\frac{d +Q_{2k-1}}{1+d
Q_{2k-1}},\,\,\,\, Q_{2k-1}=e^{2i\lambda b_1(1-x_2)}\frac{-d
+Q_{2k-2}}{1-d
Q_{2k-2}}\label{Q2k}\\
&Q_2=e^{2i\lambda b_2x_2}d.\nonumber
\end{align}

Note that if condition (\ref{cond}) is satisfied:
$b_2x_2=b_1(1-x_2),$
then  $Q_m(-i\lambda)$ as function of $\Re\lambda$ is periodic with
the period $\pi/b_2x_2,$ which implies the same property for any
solutions of equation (\ref{maineq}). The resonance spectrum is then
periodic as stated in the Introduction, Theorem \ref{main}.

\section{The poles of analytic continuation of the reflection coefficient $r_k$ to
$\cz_-$}\label{s-relation} Let $r_k$ be reflection coefficient
(\ref{rk}). In this section we express $r_k$ using the iteration
formulas in section \ref{s-banica} and extend $r_k(\lambda)$ to
$\lambda\in\cz_-.$ We show explicitly that the resonances defined in
(\ref{detres}) or (\ref{maineq}) are the poles of $r_k,$
\begin{equation*}\frac{1}{r_k}=0\,\,\Leftrightarrow\,\, \det D_{2k+1}(-i\lambda)=0
\,\,\Leftrightarrow\,\, Q_{2k}(-i\lambda)=\frac{1}{d},
\end{equation*}

We consider  the equation $P_1u=-\partial_x a(x)\partial_x
u=\lambda^2 u$ corresponding to the system with one unit-cell,
$a(x)=a_0 (x),$ for $x\in [0,1)$ and $a=1/b_1^2$ outside $[0, 1).$
Taking solution
\begin{equation}\label{koeef} u(x)=A_0 e^{i\lambda b_1x} +A_0'
e^{-i\lambda b_1 x},\,\,x<0\,\,\mbox{and}\,\,u(x)=\tilde{A}_0
e^{i\lambda b_1x} +\tilde{A}_0' e^{-i\lambda b_1 x},\,\,
x>1,\end{equation} the matching conditions imply
\begin{equation}\label{T}\left(\begin{array}{c}
    \tilde{A}_0 \\
    \tilde{A}_0'
  \end{array}\right)=T\left(\begin{array}{c}
    A_0 \\
    A_0'
  \end{array}\right),\end{equation} where $T$ is called transmission
  matrix (see \cite{Cohen-Tannoudjietal1977}).

Then coefficients of the solution of the problem with $k-$unit cells
 $P_ku=-\partial_x a(x)\partial_x
u=\lambda^2 u$ can be calculated by iteration.  Taking solution
\begin{equation}\label{coeef} u(x)=A_{-1} e^{i\lambda b_1x} +A_{-1}'
e^{-i\lambda b_1 x},\,\,x<0\,\,\mbox{and}\,\,u(x)=\tilde{A}_k
e^{i\lambda b_1x} +\tilde{A}_k' e^{-i\lambda b_1 x},\,\,
x>k,\end{equation} the coefficients of the two external regions
$x<0$ and $x>k$ are related by
\begin{equation}\label{288}\left(
    \begin{array}{c}
      A_{k} \\
      A_{k}'\\
    \end{array}
  \right)=\tilde{T}(\lambda)\left(
    \begin{array}{c}
      A_{-1} \\
      A_{-1}'\\
    \end{array}
  \right),\end{equation}
  where \begin{align}&\tilde{T}(\lambda)=\left(
                                           \begin{array}{cc}
                                             e^{-il(k+1)\lambda b_1} & 0 \\
                                             0 & e^{il(k+1)\lambda b_1} \\
                                           \end{array}
                                         \right) Q^k\left(
                                           \begin{array}{cc}
                                             e^{il\lambda b_1} & 0 \\
                                             0 & e^{-il\lambda b_1} \\
                                           \end{array}
                                         \right)=\label{299}
                                         ,\end{align}
                                         where $$Q=\left(
    \begin{array}{cc}
      Q_{11} (\lambda) & Q_{12}(\lambda) \\
      Q_{21} (\lambda) & Q_{22} (\lambda) \\
    \end{array}
  \right),\,\,\det Q=1,\,\,Q_{12}(\lambda)=Q_{21}(-\lambda),\,\,Q_{11}(-\lambda)=Q_{22}(\lambda)$$ is called
                                         iteration matrix
\begin{equation}\label{300}Q(\lambda)=D(\lambda) T(\lambda),\,\,D=\left(
                      \begin{array}{cc}
                        e^{i\lambda b_1} & 0 \\
                        0 & e^{-i\lambda b_1} \\
                      \end{array}
                    \right).
\end{equation}

Next we relate  the iteration matrix $Q$ and the monodromy matrix
$M_\lambda.$

 Relation  between the Cauchy data and coefficients in
(\ref{coeef}) is given by $$\left(
                       \begin{array}{c}
                         u(0) \\
                         \frac{(a\partial_x u)(\lambda,0)}{\lambda} \\
                       \end{array}
                     \right) =\left(
                                \begin{array}{cc}
                                  1 & 1 \\
                                  \frac{i}{b_1} & -\frac{i}{b_1} \\
                                \end{array}
                              \right)\left(
                                       \begin{array}{c}
                                         A_0 \\
                                         A_0' \\
                                       \end{array}
                                     \right)
                             =L\left(
                                       \begin{array}{c}
                                         A_0 \\
                                         A_0' \\
                                       \end{array}
                                     \right) $$ and $$\left(
                       \begin{array}{c}
                         u(1) \\
                         \frac{(a\partial_x u)(\lambda,1)}{\lambda} \\
                       \end{array}
                     \right)=\left(
                                \begin{array}{cc}
                                  1 & 1 \\
                                  \frac{i}{b_1} & -\frac{i}{b_1} \\
                                \end{array}
                              \right)\left(
                                       \begin{array}{cc}
                                         e^{i\lambda b_1 } & 0 \\
                                         0 & e^{-i\lambda b_1 } \\
                                       \end{array}
                                     \right)
                              \left(
                                       \begin{array}{c}
                                         \tilde{A}_0 \\
                                         \tilde{A}_0' \\
                                       \end{array}
                                     \right)=LD\left(
                                       \begin{array}{c}
                                         \tilde{A}_0 \\
                                         \tilde{A}_0' \\
                                       \end{array}
                                     \right).$$

Using
$$M_\lambda\left(
                       \begin{array}{c}
                         u(0) \\
                         \frac{(a\partial_x u)(\lambda,0)}{\lambda} \\
                       \end{array}
                     \right)\,\,=\,\,\left(
                       \begin{array}{c}
                         u(1) \\
                         \frac{(a\partial_x u)(\lambda,1)}{\lambda} \\
                       \end{array}
                     \right)$$ and definition of $T,$ equation (\ref{T}), we get
                     \begin{equation}\label{similar}M_\lambda=LDTL^{-1}=LQL^{-1},\,\,\mbox{with}\,\,L=\left(
      \begin{array}{cc}
        1 & 1 \\
        \frac{i}{b_1} & -\frac{i}{b_1} \\
      \end{array}
    \right).\end{equation}

We have
\begin{align*}
&Q_{12}(\lambda)=e^{i\lambda (1-x_2) b_1}\frac{b_2^2-b_1^2}{4b_1
b_2}\left(e^{-i\lambda b_2 x_2} -e^{i\lambda b_2
x_2}\right),\\
&Q_{22}(\lambda)=e^{-i\lambda b_1 (1-x_2)}\left( (r+1)e^{-i\lambda
b_2 x_2} -(r-1) e^{i\lambda b_2 x_2}\right) .
\end{align*}

As \begin{align*} &\det \tilde{D}_3 =e^{-i\lambda b_1
x_2}(b_2^2-b_1^2)(e^{-i\lambda b_2 x_2} - e^{i\lambda b_2
x_2}),\\
&\det D_3=e^{i\lambda b_1 x_2}\left( (b_2-b_1)^2 e^{i\lambda b_2
x_2}- (b_2 +b_1)^2 e^{-i\lambda b_2 x_2}\right),
\end{align*}
then
\begin{equation}\label{q1222}
Q_{12}(\lambda)= \frac{e^{i\lambda b_1}}{4b_1b_2}\det
\tilde{D}_3,\,\,Q_{22}(\lambda)=-\frac{e^{-i\lambda
b_1}}{4b_1b_2}\det D_3.\end{equation}

Using the iteration relations between $\det
\tilde{D}_{2k-1},\,\,\det D_{2k-1}$ and $\det
\tilde{D}_{2k+1},\,\,\det D_{2k+1}$ we get
\begin{lemma}
We have $$ Q^{k}=\left(
    \begin{array}{cc}
      Q_{22}^{(k)} (-\lambda) & Q_{12}^{(k)}(\lambda) \\
      Q_{12}^{(k)} (-\lambda) & Q_{22}^{(k)}(\lambda) \\
    \end{array}
  \right),$$ where
$$ Q_{12}^{(k)}=\frac{e^{ki\lambda b_1}}{(4b_1b_2)^k}\det
\tilde{D}_{2k+1},\,\,Q_{22}^{(k)} =-\frac{e^{-ik\lambda
b_1}}{(4b_1b_2)^k}\det D_{2k+1}.$$
\end{lemma}

Using (\ref{similar}) we have also
$$Q^k=L^{-1}M_\lambda^k L=\left(
                 \begin{array}{cc}
                   \frac12 (\alpha_k +\beta_k\frac{i}{b_1})-\frac{ib_1}{2} (\gamma_k +\delta_k\frac{i}{b_1}) & \frac12 (\alpha_k -\beta_k\frac{i}{b_1})-\frac{ib_1}{2}
                   (\gamma_k
                   -\delta_k\frac{i}{b_1})\\ \\
                   \frac12 (\alpha_k +\beta_k\frac{i}{b_1})+\frac{ib_1}{2} (\gamma_k +\delta_k\frac{i}{b_1}) & \frac12 (\alpha_k -\beta_k\frac{i}{b_1})+\frac{ib_1}{2} (\gamma_k -\delta_k\frac{i}{b_1}) \\
                 \end{array}
               \right).$$

Using equation (\ref{rk}) we get
$$r_k=-\frac{(\alpha_k -\delta_k) +i (b_1\gamma_k
+\frac{\beta_k}{b_1})}{\alpha_k +\delta_k +i(b_1\gamma_k
-\frac{\beta_k}{b_1})}=-\frac{2Q_{21}^{(k)}(\lambda)}{2Q_{22}^{(k)}(\lambda)}
= -\frac{Q_{21}^{(k)}(\lambda)}{Q_{22}^{(k)}(\lambda)}.$$

Let first $k=1.$ Using equation (\ref{q1222}) we get  $$ Q_{21}
(\lambda)=-e^{-i2\lambda (1-x_2)b_1}
                    Q_{12}=-\frac{e^{-i\lambda (1-x_2)b_1}e^{i\lambda b_1
x_2}}{4b_1b_2} \det \tilde{D}_3 (\lambda)$$

 Hence
\begin{align*}r_1&=-\frac{e^{-2i\lambda b_1 (1-x_2)}e^{i\lambda b_1(1-x_2)}
e^{i\lambda b_1 x_2}\det \tilde{D}_3}{e^{-i\lambda b_1}\det D_3}
=\\
&=-e^{-2i\lambda b_1 (1-x_2)} e^{2i\lambda b_1 x_2}\frac{\det
\tilde{D}_3}{\det D_3}= -e^{-2i\lambda b_1 (1-x_2)}Q_3
(-i\lambda).\end{align*}

For $k>1$ we use $$Q^k=\frac{\sin k\theta l}{\sin \theta l} Q
-\frac{\sin (k-1)\theta l}{\sin\theta l} I=\left(
                     \begin{array}{cc}
                       \frac{\sin k\theta l}{\sin \theta l} Q_{11} -\frac{\sin (k-1)\theta
l}{\sin\theta l}  &  \frac{\sin k\theta l}{\sin \theta l} Q_{12}
 \\
                       \frac{\sin k\theta l}{\sin \theta l} Q_{21} &  \frac{\sin k\theta l}{\sin \theta l} Q_{22} -\frac{\sin (k-1)\theta
l}{\sin\theta l} \\
                     \end{array}
                   \right)$$
                    and get
                    $$Q_{21}^{(k)}=-e^{-i2\lambda (1-x_2)b_1}
                    Q_{12}^{(k)}=-e^{-i2\lambda
                    (1-x_2)b_1}\frac{e^{ki\lambda b_1}}{(4b_1b_2)^k}
                    \det \tilde{D}_{2k+1}.$$
                    Hence, by (\ref{rk}),
                    $$r_k=-\frac{Q_{21}^{(k)}}{Q_{22}^{(k)}}=-\frac{e^{-i2\lambda
                    (1-x_2)b_1}e^{ki\lambda b_1}
                    \det \tilde{D}_{2k+1}}{e^{-ik\lambda b_1}\det
                    D_{2k+1}}=-e^{-i2\lambda (1-x_x)b_1} e^{2ki\lambda b_1}\frac{\det \tilde{D}_{2k+1}}{\det D_{2k+1}}.$$

We have
\begin{equation}\label{rk_linear_fractional}r_k=-e^{-2i\lambda b_1
(1-x_2)}Q_{2k+1} (-i\lambda) =-\frac{-d +Q_{2k} (-i\lambda)}{1-d
Q_{2k} (-i\lambda)}.\end{equation} Thus we have proved that the
poles of $r_k$ on $\cz_-$ coincide with the solutions of
(\ref{maineq}):
$$1-d Q_{2k} (-i\lambda)=0\,\,\Leftrightarrow\,\, Q_{2k} (-i\lambda)
=\frac{1}{d}.$$ Hence, we have
\begin{proposition}\label{prop_rk}
The reflection coefficient $r_k(\lambda)$ continuous to $\cz_-$ with
the poles at the resonance spectrum ${\rm Res}\,(P_k).$
\end{proposition}

\section{Convergence of the resonances to the real axis as
$k\rightarrow\infty$}\label{s-convergence} In this section we prove
that the resonances spectrum ${\rm Res}(P_k)$ converges to the real
axes as $k\rightarrow\infty.$

First we note that   the function $\theta (\lambda)=\arccos\left(
(\alpha(\lambda) +\delta (\lambda))/2\right)$ has analytic
continuation onto the domain $\cz\setminus \bigcup_{n=1}^\infty g_n$
by the formula $\theta(\overline{\lambda})=\overline{\theta
(\lambda)}$
and $\Im \theta (\lambda) <0$ for $\Im\lambda <0.$
 For
$\lambda\in\bigcup_{n=1}^\infty g_n$ we set
$\theta(\lambda)=\theta(\lambda-i0)$ and $\theta$ is pure imaginary
there.

By Proposition \ref{prop_rk} the function $r_k(\lambda)$ is analytic
in $\lambda\in\cz_-\setminus{\rm Res}\,(P_k),$ and the formula
\begin{equation*} r_k(\lambda)=-\frac{(\alpha-\delta)+i(b_1\gamma
+\frac{\beta}{b_1} )}{2\sin l\theta (\lambda)\frac{\cos
kl\theta(\lambda)}{\sin kl\theta(\lambda)} +i(b_1\gamma
-\frac{\beta}{b_1})}
\end{equation*}
extends to  $\cz_-\setminus{\rm Res}\,(P_k),$ where we have
$\Im\theta (\lambda) <0.$

Furthermore,  if $\Im\theta (\lambda) <0,$ then
$$\frac{\cos k\theta(\lambda)}{\sin k\theta (\lambda)}\rightarrow i
\,\,\mbox{as}\,\, k\rightarrow \infty, $$ and we have
$$r_k\rightarrow \frac{i(\alpha-\delta)-(b_1\gamma
+\frac{\beta}{b_1} )}{2\sin l\theta (\lambda) +(b_1\gamma
-\frac{\beta}{b_1})}=\tilde{r}(\lambda).$$

Note that for $\lambda\in\rz,$
$\tilde{r}(\lambda)=\overline{r(-\lambda)}.$

The limit extends also to $\lambda\in \cup_{n=1}^\infty g_n,$ where
$\theta(\lambda)$ is pure imaginary.

Let $g_k(\lambda),$ $g(\lambda)$ denote the denominators of $r_k$
respectively $\tilde{r}$ for $\lambda\in\cz_-.$ Then
$$g_k(\lambda)\,\,\rightarrow\,\, g(\lambda)=2\sin l\theta (\lambda) +\left(
b_1\gamma -\frac{\beta}{b_1}\right)$$ uniformly on any compact
subsets of $\cz_-$ and the limiting function $g$ is analytic
 on
$\cz_-.$

For any $k$ all zeros of $g_k$
 have
negative imaginary part: ${\rm Res}\,(P_k)\in\cz_-.$ Then, by the
Hurwitz's theorem, the zeros of $g_k$ can only converge to the real
axis.
\hfill\qed

\section{The limit  $\lim_{k\rightarrow\infty}r_k$ as a fixed point of a sequence of linear-fractional automorphisms of the unit
disk}\label{s-fixed_point} Equation (\ref{rk_linear_fractional})
together with iteration relations (\ref{rk_linear_fractional}) shows
that for $\lambda\in\rz$ real and  for all $k=1,2,3,\ldots$ the
reflection coefficients for the $k+1$ and $k$ cells media are
related by $r_{k+1}=f_\lambda (r_k),$ where $f_\lambda$ is
linear-fraction automorphism of the unit disk. Hence we get
\begin{equation}\label{iterates}r_k=f_\lambda^{[k]}(r_1):=\underbrace{f_\lambda\circ f_\lambda\ldots
f_\lambda}_{k\,\,\mbox{iterates}} (r_1).\end{equation} Here $r_1$ is
the reflection coefficient for $P_1,$
$$r_1=-\frac{(1-e^{2i\lambda b_2 x_2})d}{1-d^2 e^{2i\lambda b_2
x_2}}$$ and for $\lambda\in\rz$ we have $$|r_1|^2\leq
\frac{4d^2}{1+4d^2+d^4} <1.$$ In this section we will consider the
limit of the sequence $f_\lambda^{[k]}$ as $k\rightarrow\infty.$

First we recall some well-known facts   on the convergence behavior
of a sequence $\{ f^{[k]}\}$ when $f$ is general linear-fractional
automorphism $f$ of the unit disc $D=\{ z\in\cz;\,\,|z|<1\},$
$$f(z)=\frac{b -z}{1-\overline{b}z},\,\,b\in
D,\,\,z\in\cz\setminus\{1/\overline{b}\}.$$ We refer to the  paper
of Burckel \cite{Burckel1981}) for the details.

 In general
situation ($b\neq 0$), $f$ has two fixed points $z_1,$ $z_2.$ There
are three cases to consider.
\begin{description}
\item[Hyperbolic] $f$ has two (distinct) fixed points on $\partial D$
 $$|z_1|=|z_2|=1,\,\,z_1\neq
z_2.$$ In this case the sequence $\{f^{[k]}\}$ converges uniformly
on compact subsets in $D$ to one of these points.
\item[Parabolic] $f$ has  one (double) fixed point on $\partial D$
 $$z_1=z_2,\,\, \mbox{with}\,\,|z_1|=1.$$ In this case the sequence
$\{f^{[k]}\}$ converges uniformly on compact subsets in $D$ to this
fixed point.
\item[Elliptic] $f$ has two fixed points: one fixed point $z_1\in D$ and one fixed point $z_2=1/\overline{z}\not\in D.$
In this case either $f$ is periodic in the sense that $f^{[n]}=I$
for some $n,$ or the orbit $\{f^{[k]};\,\,n\in\nz\}$ is dense in the
compact group of all conformal automorphisms of $D$ which fix $z_1.$
\end{description}

We apply these results  to $f=f_\lambda.$  In order to simplify the
formulas we suppose
\begin{equation}\label{co}b_2 x_2=b_1(1-x_2).\end{equation}

Then we have the following expression for $f_\lambda:$
$$
f_\lambda(z)=\frac{-d+\eta\frac{d+\eta z}{1+d\eta z}}{1-d\eta
\frac{d+\eta z}{1+d\eta z}},\,\,\mbox{where}\,\,\eta:=e^{2i\lambda
b_2 x_2}.$$

In general situation $\eta\neq 1,$ equation $f_\lambda(z)=z$ have
solutions
\begin{equation}\label{fixedpoints}
z_{1,2}=\frac{-(1+\eta)\pm\sqrt{(1+\eta)^2-4d^2\eta}}{2d\eta}=\frac{-\cos
(\lambda b_2 x_2)\pm\sqrt{\cos^2(\lambda b_2 x_2)-d^2}}{d
e^{i\lambda b_2 x_2}}.
\end{equation}
Note that
\begin{align}
\mbox{if}\,\,& d^2<\cos^2(\lambda b_2 x_2) <1\,\,\mbox{then}\,\,
|z_1|<1,\,\,|z_2|>1\,\,\mbox{and}\,\,z_1\cdot\overline{z_2}=1\label{eli}\\
\mbox{if}\,\,& \cos^2(\lambda b_2 x_2) \leq d^2\,\,\mbox{then}\,\,
|z_1|=|z_2|=1,\label{hyppar}
\end{align}
and $z_1=z_2$ if $\cos^2(\lambda b_2 x_2)=d.$

 Note that if (\ref{co}) is satisfied then the
Lyapunov function is given by (\ref{Lyp}),
$F(\lambda)=(\rho+1)\cos^2(\lambda b_2x_2) -\rho.$ The inner points
of the allowed bands $\bigcup_{n=1}^\infty b_n$ satisfy  $
|F(\lambda)|^2 < 1$ which is equivalent to $d^2<\cos^2(\lambda b_2
x_2) < 1.$ On the spectral gaps $\bigcup_{n=1}^\infty g_n$ we have
$\cos^2(\lambda b_2 x_2)<d^2 .$ The  non-generated edge point of a
band is given by $\cos^2(\lambda b_2 x_2)=d^2.$ The degenerated edge
point $\lambda_0$ satisfy $\cos^2(\lambda b_2
x_2)=1\,\,\Leftrightarrow\,\,\lambda_0=\pi  m/b_2x_2.$

Hence we get the following version of the forth result in Theorem
\ref{tMV}:
\begin{proposition} If $\lambda\in\bigcup_{n=1}^\infty g_n$ then
the linear fractional automorphism $f_\lambda$ has two fixed points
$z_1,$ $z_2$ of hyperbolic type: $|z_1|=|z_2|=1.$ The sequence $r_k$
converges to $r$  as $k\rightarrow\infty,$ where
$$r=\lim_{k\rightarrow\infty} f_\lambda^{[k]}(r_1)$$
 is either $z_1$ or  $z_2.$

 If $\lambda$ is a non-degenerated band edge point: $ F(\lambda)=\pm 1$ and $F'(\lambda)\neq 0,$ then $f_\lambda$ has one (double) fixed
 point $z_1=z_2$ and $|z_1|=1.$

 We have $$r=z_1=\lim_{k\rightarrow\infty} r_k=\lim_{k\rightarrow\infty} f_\lambda^{[k]}(r_1).$$

If $\lambda=\lambda_0$ is degenerate band edge point $
F(\lambda_0)=\pm 1,$ $F'(\lambda_0)= 0$ (degenerate band edge) then
$r_k(\lambda)=0$ for all $k=1,2,3\ldots.$
\end{proposition}

The Proposition is still valid if condition (\ref{co}) is not
imposed.


\appendix
\section{Numerical calculations}\label{s-numerics}
In this section we present some examples of the resonance spectrum.
The resonances are solutions of the equations
$D_{2k+1}(-i\lambda)=0,$ where $D_n,$ $n=2k+1,$ are defined
iteratively by (\ref{indrel}), Section \ref{s-banica}. The zeros of
$D_{2k+1}(-i\lambda)$ are calculated numerically by using the Newton
procedure.
  Using Matlab we plot the resonance spectrum for
the number of identical cells $k=3,4,5.$  In the same figure we show
the band spectrum for the corresponding periodic problem satisfying
(\ref{bandspectrum}). The small circles on the real axis marks the
position of $\lambda_0=\pi m/ (x_2 b_2)$ when the one-cell system is
perfectly transparent: $|t_1(\lambda_0)|^2=1,$
$|r_1(\lambda_0|^2=0.$

\begin{figure}[htbp]
\includegraphics*[width =13 cm]{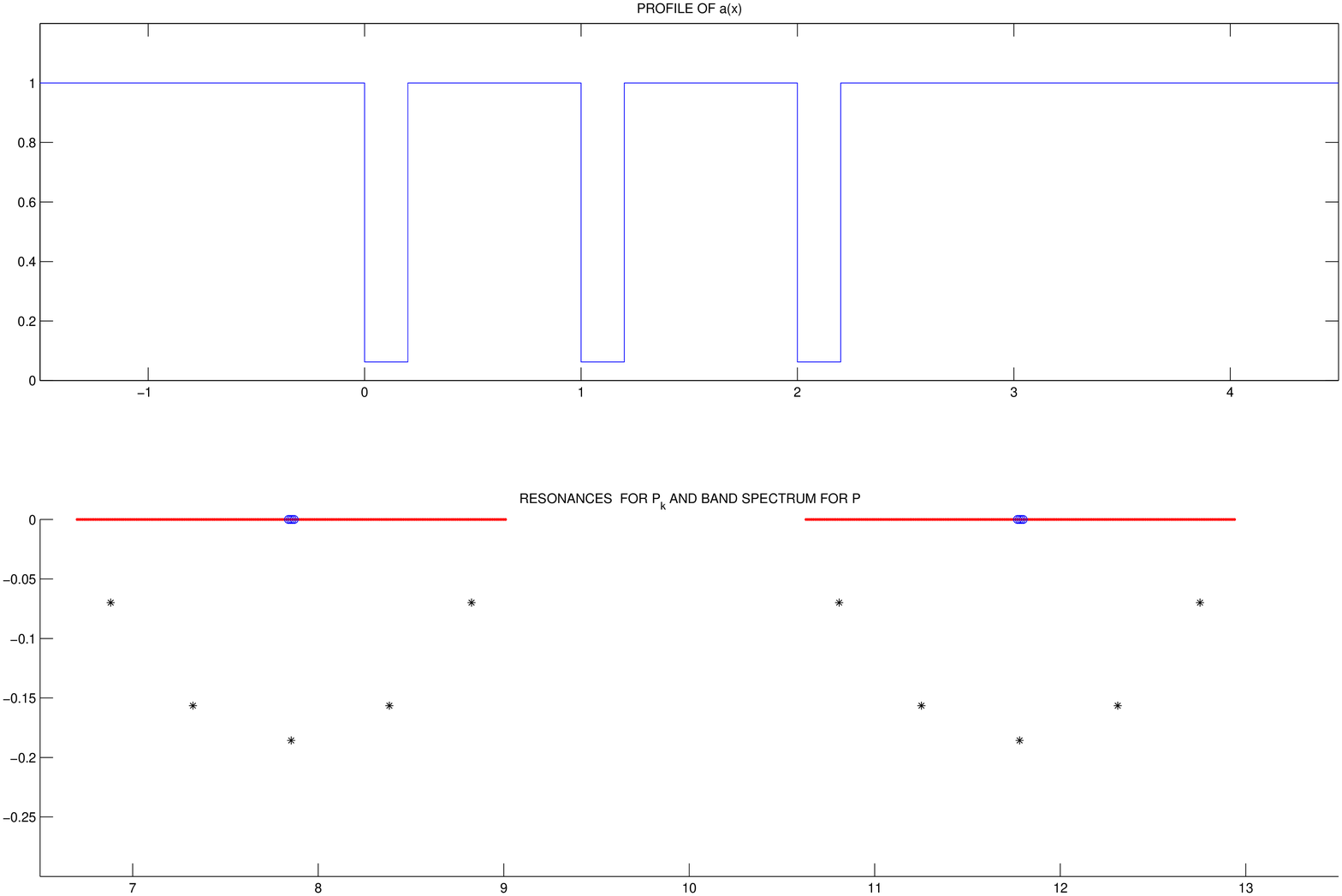}

\includegraphics*[width =13 cm]{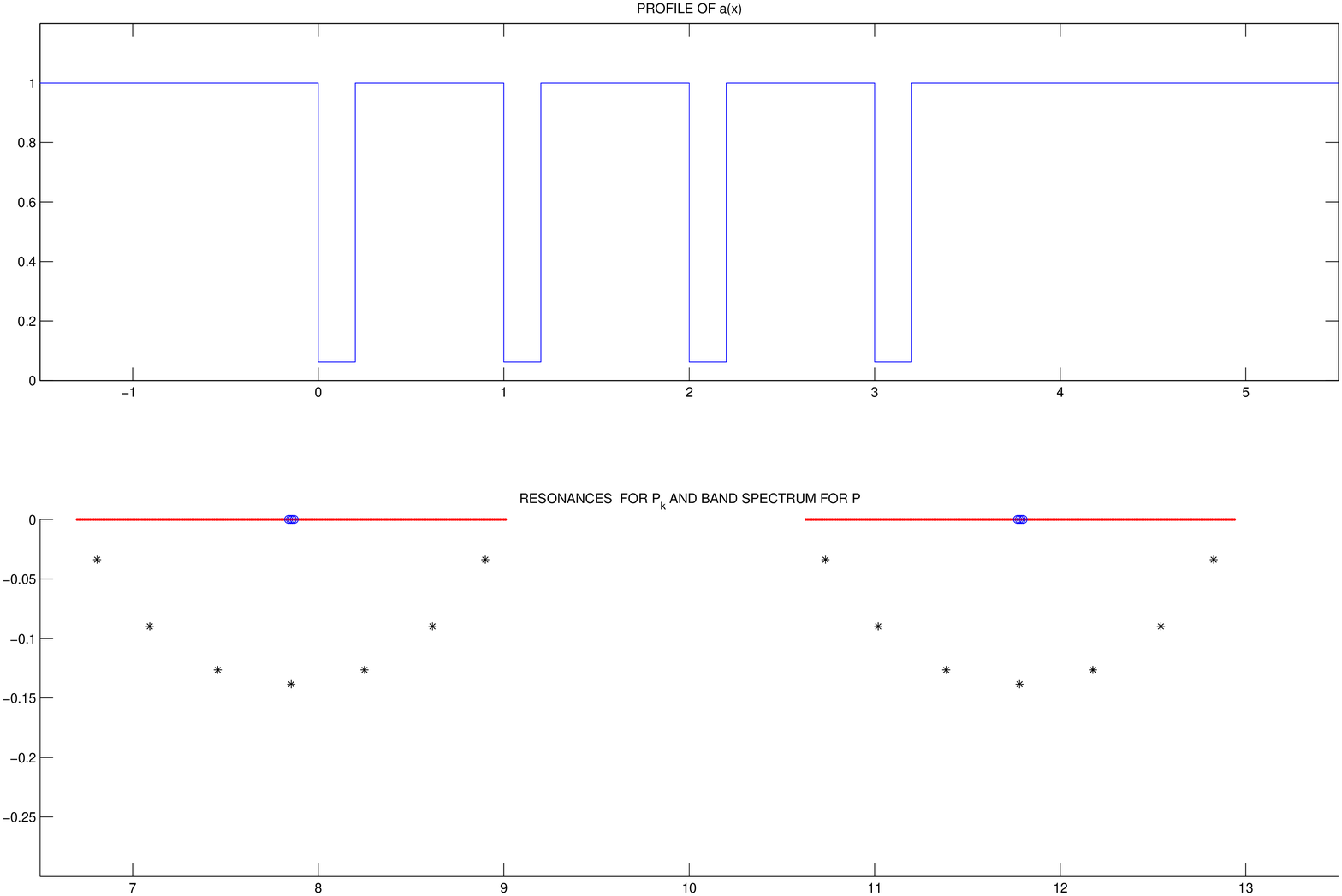}

\includegraphics*[width =13 cm]{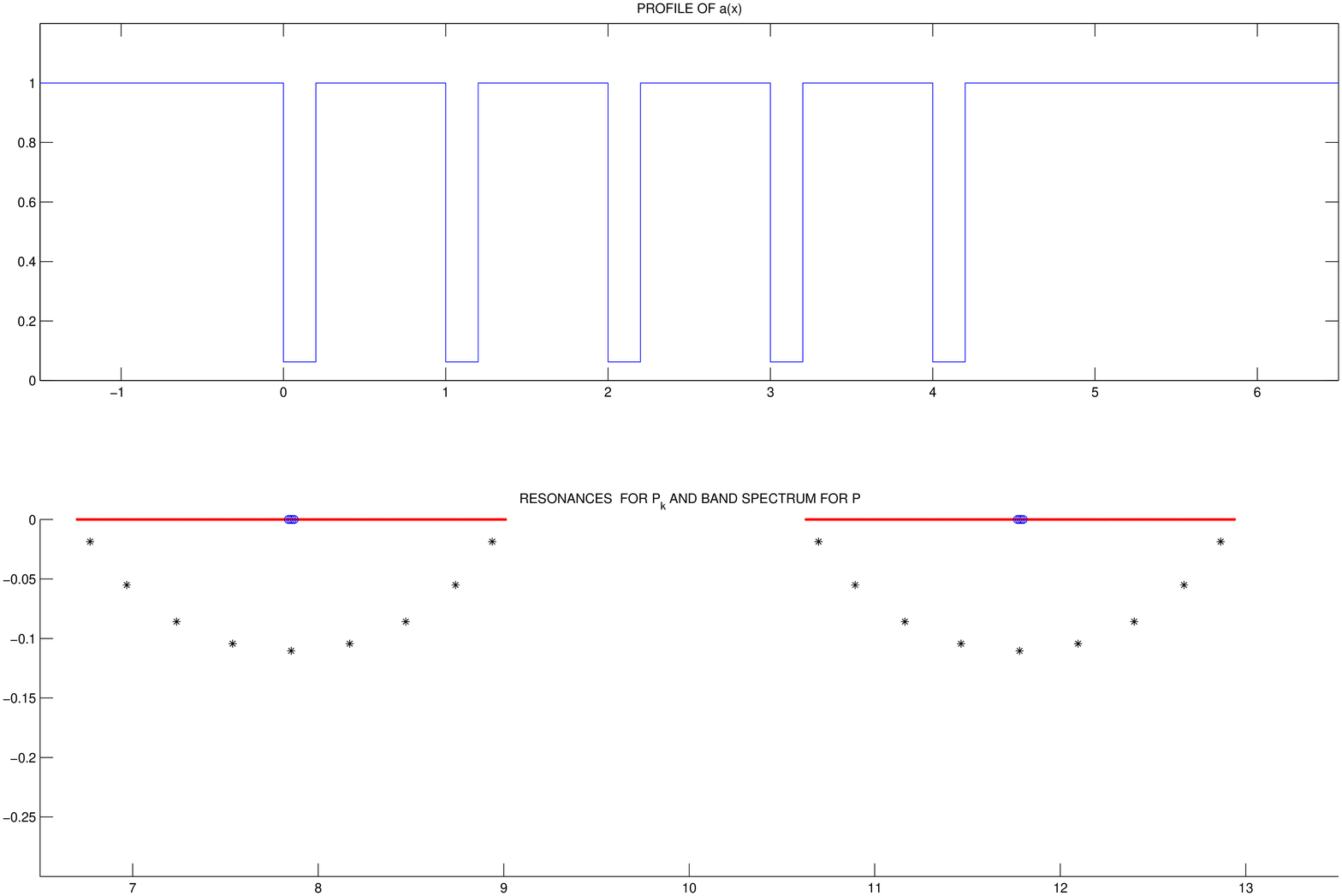}
\caption{\label{Fig1}{\em $b_1=1,\,\,b_2=4,\,\,x_2=0.2$ } }
\end{figure}

On Figure (\ref{Fig1}) condition (\ref{cond}) is satisfied:
$b_1=1,\,\,b_2=4,\,\,x_2=0.2.$

\begin{figure}[htbp]
\includegraphics*[width =13 cm]{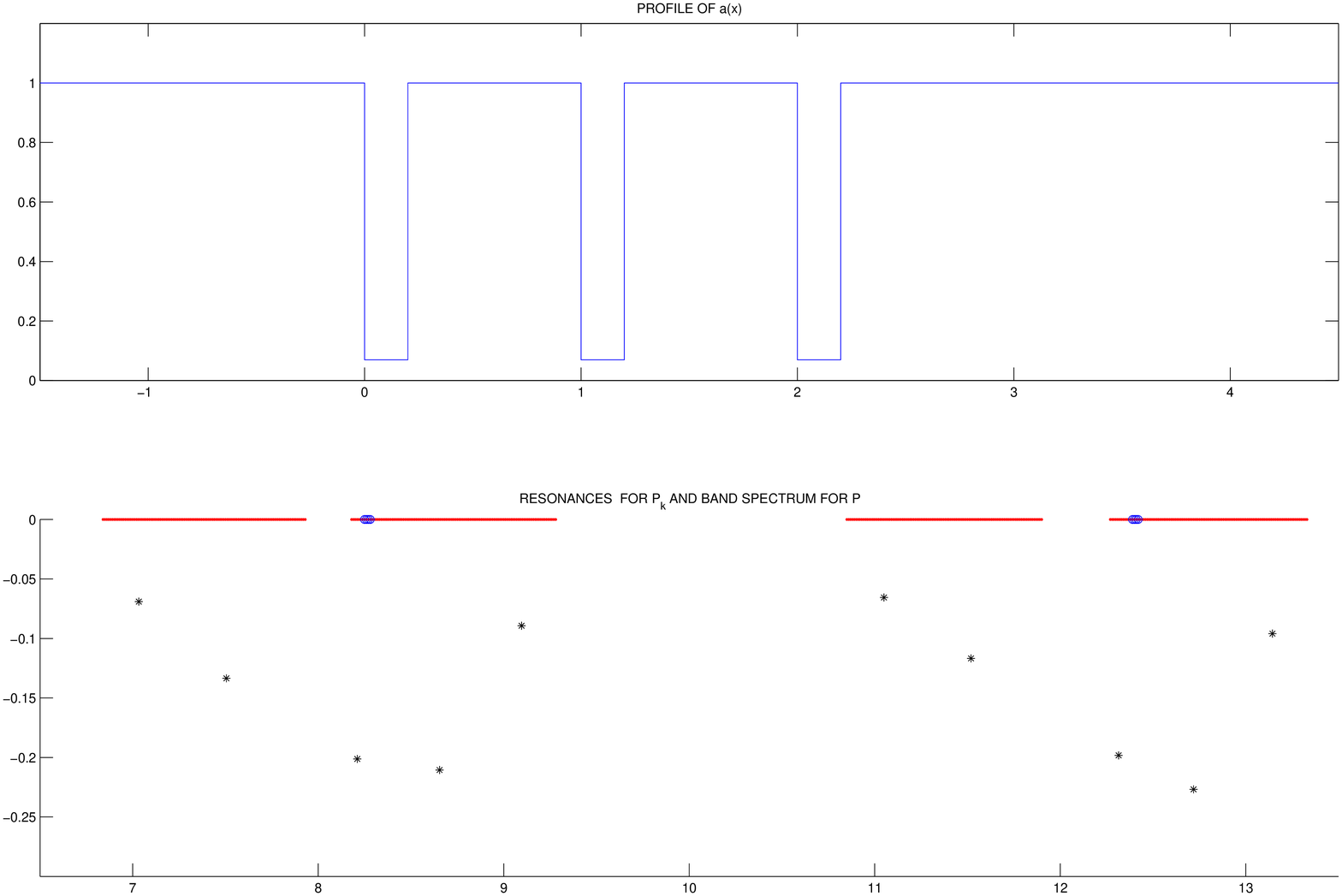}

\includegraphics*[width =13 cm]{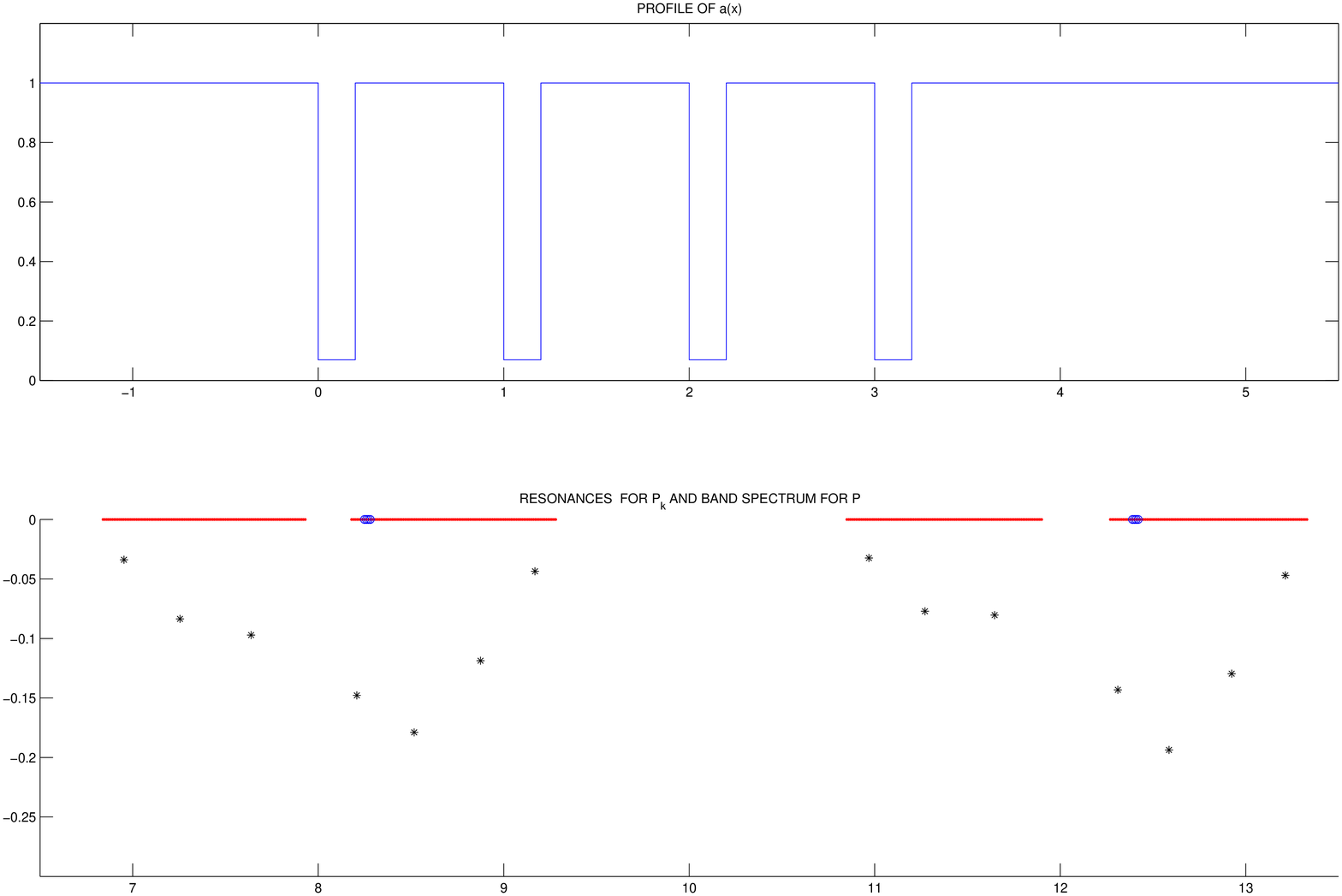}

\includegraphics*[width =13 cm]{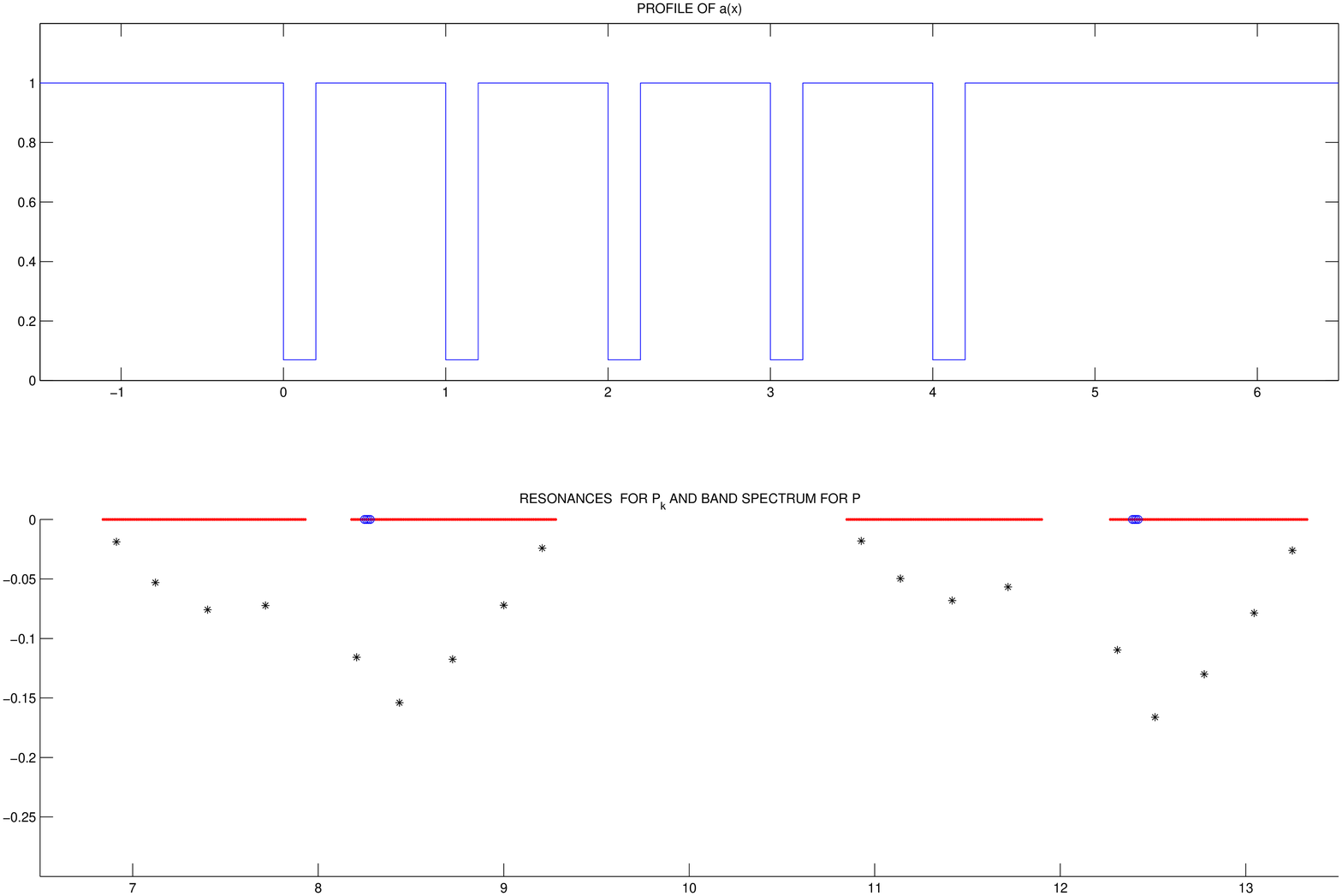}
\caption{\label{Fig2}{\em $b_1=1,\,\,b_2=3.8,\,\,x_2=0.2$} }
\end{figure}

On Figure (\ref{Fig2}) condition (\ref{cond}) is not satisfied:
 $b_1=1,\,\,b_2=3.8,\,\,x_2=0.2.$

\begin{figure}[htbp]
\includegraphics*[width =13 cm]{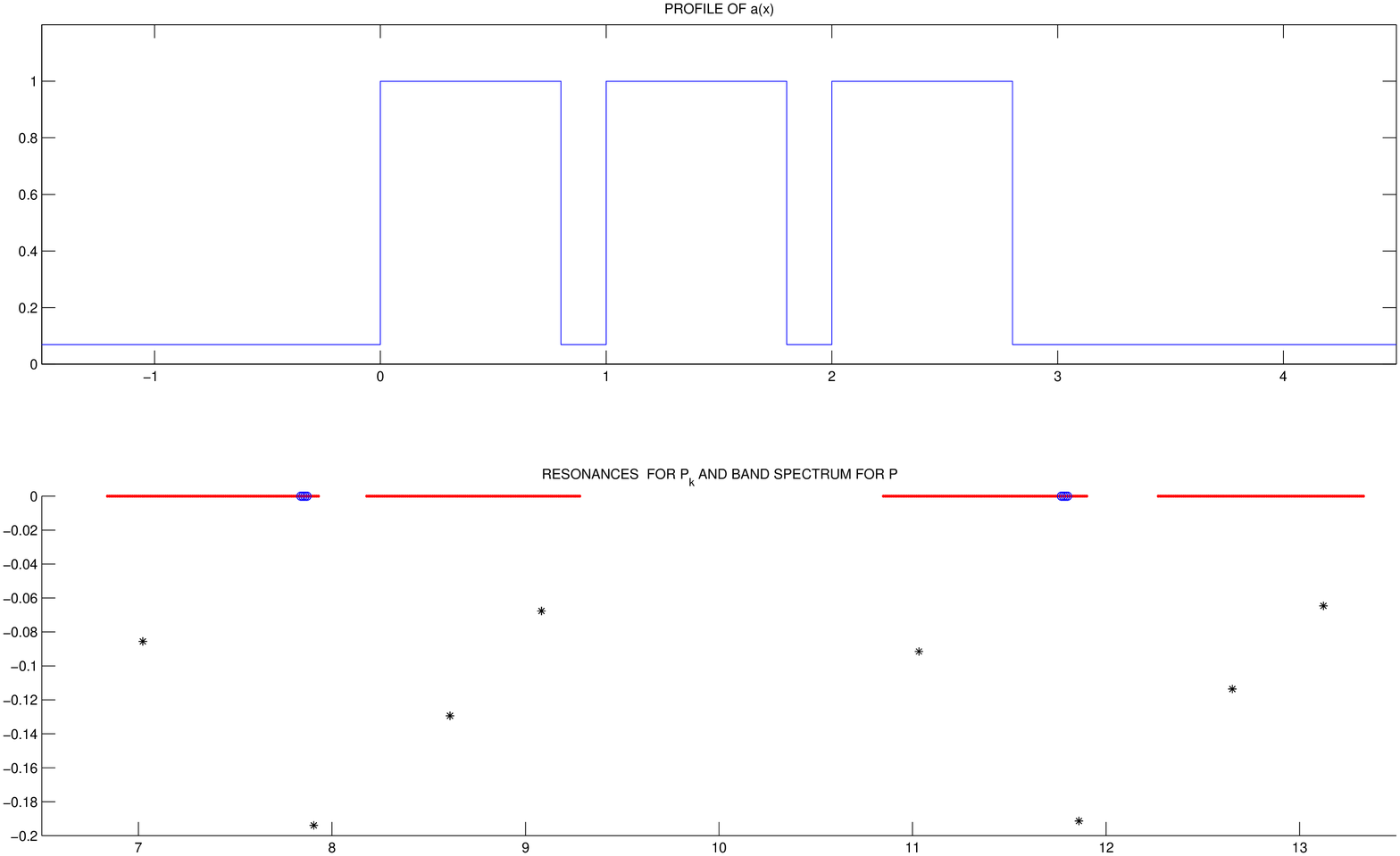}

\includegraphics*[width =13 cm]{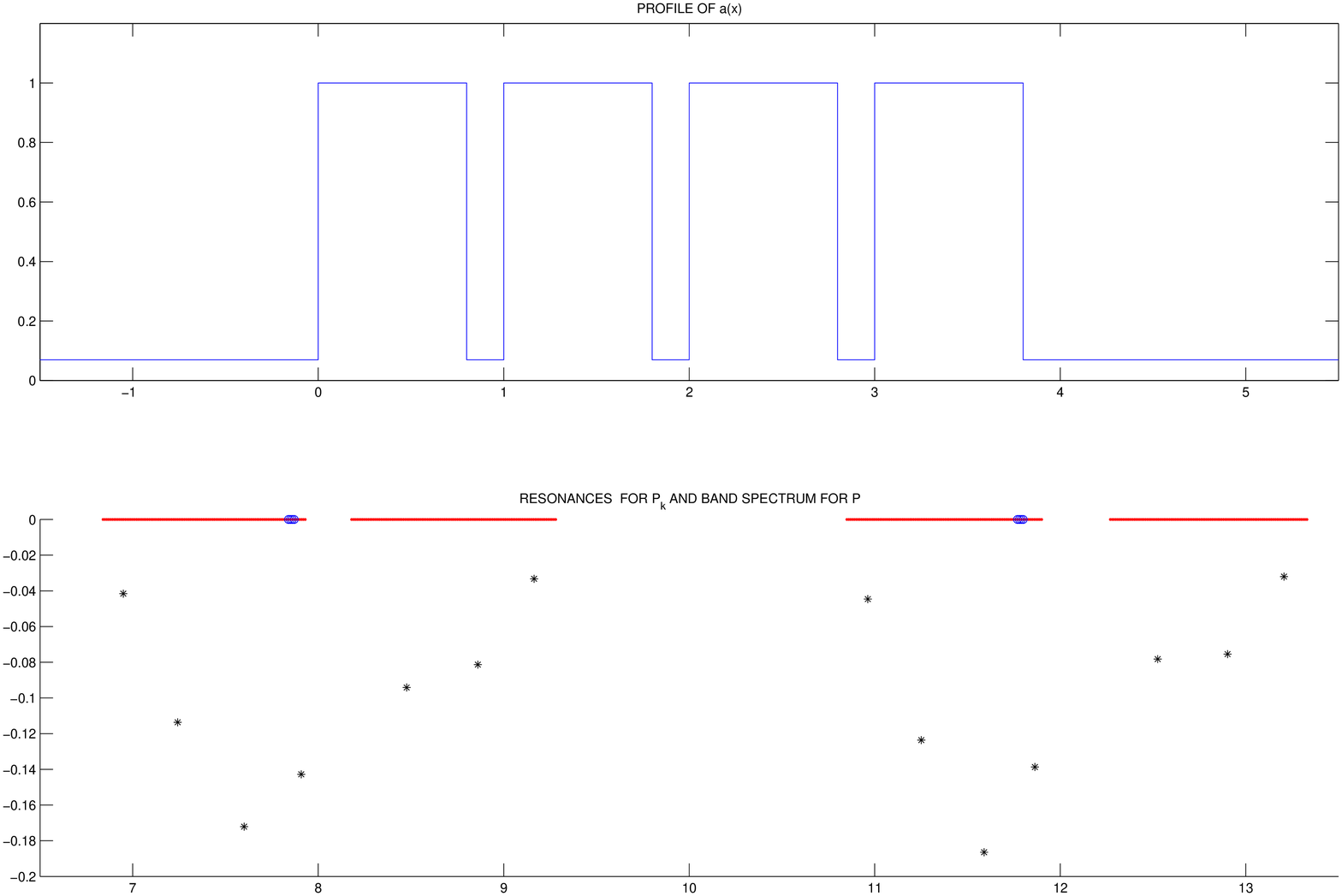}

\includegraphics*[width =13 cm]{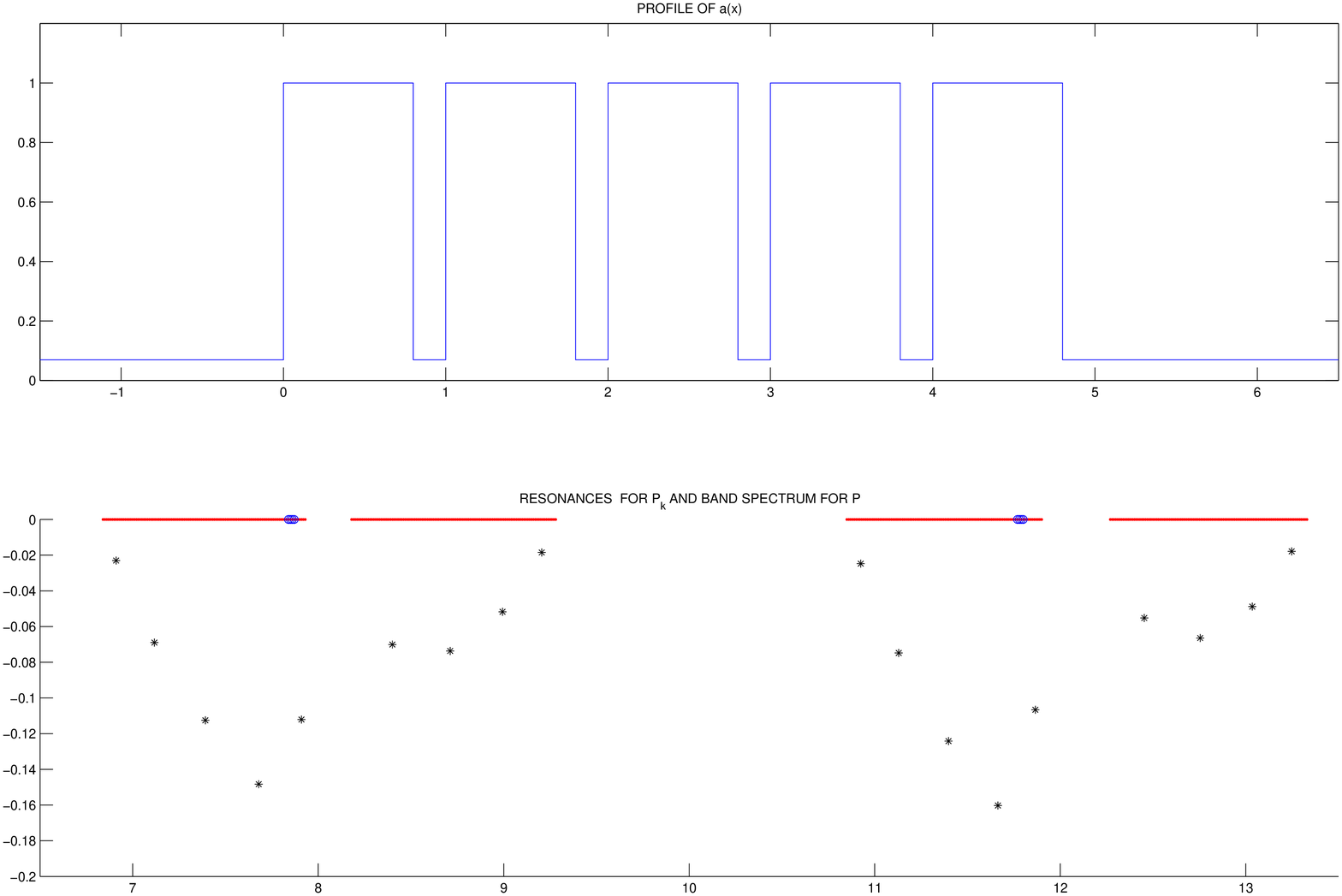}
\caption{\label{Fig3}{\em $b_1=3.8,\,\,b_2=1,\,\,x_2=0.8$} }
\end{figure}

On Figure (\ref{Fig3}) condition (\ref{cond}) is not satisfied:
$b_1=3.8,\,\,b_2=1,\,\,x_2=0.8.$




\end{document}